\documentclass[twoside,12pt]{article}
\usepackage{geometry}
\usepackage{amsmath,amssymb,amsfonts, amsthm}
\usepackage{txfonts}
\usepackage{todonotes}
\usepackage{graphicx}
\usepackage{mathrsfs}
\usepackage{color}
\usepackage{float}
\usepackage{extarrows}
\usepackage{exscale}
\usepackage{titlesec}
\titleformat{\section}[block]{\large\center\sc}{\arabic{section}}{0.5em}{}[] 
\usepackage{relsize}
\usepackage{color}
\usepackage{hyperref}
\usepackage[titles]{tocloft} 
\usepackage{fancyhdr}
\theoremstyle{plain}
\newtheorem{theorem}{Theorem}[section]
\newtheorem{lemma}[theorem]{Lemma}

\newtheorem{conjecture}[theorem]{Conjecture}

\theoremstyle{definition}

\theoremstyle{remark}
\newtheorem{remark}[theorem]{Remark}

\let\oldsection\section
\renewcommand\section{\setcounter{equation}{0}\oldsection}

\def\be{\begin{equation}}
\def\ee{\end{equation}}
\def\bes{\begin{equation*}}
\def\ees{\end{equation*}}
\def\bs{\begin{split}}
\def\es{\end{split}}
\def\bali{\begin{aligned}}
\def\eali{\end{aligned}}
\newcommand{\pf}{\noindent {\bf Proof. \hspace{2mm}}}
\def\bR{{\mathbb R}}
\def\un{\underbrace}
\def\al{\alpha}

\def\ve{\varepsilon}

\def\t{\tilde}
\def\th{\theta}

\def\g{\gamma}

\def\dl{\delta}
\def\Dl{\Delta}
\def\lt{\left}
\def\rt{\right}
\def\les{\lesssim}
\def\si{\sigma}

\def\i{\infty}
\def\p{\partial}
\def\f{\frac}
\def\na{\nabla}
\def\o{\omega}
\def\O{\Omega}
\def\s{\sqrt}
\def\q{\quad}

\def\bl{\boldsymbol}
\def\mO{\mathcal{O}}

\allowdisplaybreaks
\def\mD{\mathcal{D}}

\begin{document}

\title{\normalsize\bf CHARACTERIZATION OF SMOOTH SOLUTIONS TO THE NAVIER-STOKES EQUATIONS IN A PIPE WITH TWO TYPES OF SLIP BOUNDARY CONDITIONS}

\author{\normalsize\sc Zijin Li, Xinghong Pan and Jiaqi Yang}

\date{}

\maketitle

\begin{abstract}

Smooth solutions of the stationary Navier-Stokes equations in an infinitely long pipe, equipped with the Navier-slip or Navier-Hodge-Lions boundary condition, are considered in this paper. Three main results are presented.

{First, when equipped with the Navier-slip boundary condition,} it is shown that, $W^{1,\i}$ axially symmetric solutions with zero flux at one cross section, must be swirling solutions: $u=(- C x_2, C x_1,0)$, and $x_3-$periodic solutions must be helical solutions: $u=(-C_1x_2,C_1x_1,C_2)$.

Second, also equipped with the Navier-slip boundary condition, if the swirl or vertical component of the axially symmetric solution is independent of the vertical variable $x_3$, solutions are also proven to be helical solutions. In the case of the vertical component being independent of $x_3$, the $W^{1,\i}$ assumption is not needed. In the case of the swirl component being independent of $x_3$, the $W^{1,\i}$ assumption can be relaxed extensively such that the horizontal radial component of the velocity, $u_r$, can grow exponentially with respect to the distance to the origin. Also, by constructing a counterexample, we show that the growing assumption on $u_r$ is optimal.

Third, when equipped with the Navier-Hodge-Lions boundary condition, we can show that if the gradient of the velocity grows sublinearly, then the solution, enjoying the Liouville-type theorem, is a trivial shear flow: $(0,0,C)$.

\medskip

{\sc Keywords:} Navier-Stokes system, Navier-slip boundary, Navier-Hodge-Lions boundary, axially symmetric, helical solutions.

{\sc Mathematical Subject Classification 2020:} 35Q35, 76D05

\end{abstract}


\section{Introduction}
\q\ The 3D stationary Navier-Stokes (NS) equations which describes the motion of stationary viscous incompressible fluids follows that
\be\label{NS}
\lt\{
\begin{aligned}
&u\cdot\na u+\na p-\Dl u=0,\\
&\na\cdot u=0,
\end{aligned}
\rt.\q \text{in}\q \mathcal{D}\subset \bR^3.
\ee
Here  $u(x)\in\mathbb{R}^3$, $p(x)\in\mathbb{R}$ represents the velocity and the scalar pressure respectively. In this paper, we consider the domain $\mD$ to be an infinitely long pipe, i.e.
\bes
\mD=\lt\{x:\, |x_h|< 1, x_3\in \bR\rt\},
\ees
where $x=(x_1,x_2,x_3)$, $x_h=(x_1,x_2)$ and $|x_h|=\s{x^2_1+x^2_2}$. The boundary condition will be equipped with the following:

{\noindent\bf The total Navier-slip boundary condition:}
\be\label{NSBC}\tag{NSB}
\left\{
\begin{aligned}
&(\mathbb{S}u\cdot\boldsymbol{n})_\tau=0,\\
&u\cdot\boldsymbol{n}=0,\\
\end{aligned}
\right.\quad\forall x\in\p \mD,
\ee
{\noindent or {\bf the Navier-Hodge-Lions boundary condition}:

\be\label{SLIP}\tag{NHLB}
\left\{
\begin{aligned}
&(\na\times u)\times \boldsymbol{n}=0,\\
&u\cdot\boldsymbol{n}=0,\\
\end{aligned}
\right.\quad\forall x\in\p \mD,
\ee
}
Here $\mathbb{S}u=\frac{1}{2}\left(\nabla u+(\nabla u)^T\right)$ is the stress tensor, where $(\nabla u)^T$ is the transpose of the Jacobian matrix $\nabla u$, and $\boldsymbol{n}$ is the unit outer normal vector of $\p \mD$. For a vector field $v$, $v_\tau$ stands for its tangential part: $v_\tau:=v-(v\cdot\boldsymbol{n})\boldsymbol{n}.$ The condition \eqref{NSBC} is from the general Navier-slip boundary condition and impermeable boundary condition which was introduced by Claude-Louis Naiver in 1820s \cite{Navier}:
\be\label{NBC}
\left\{
\begin{aligned}
&2(\mathbb{S}u\cdot\boldsymbol{n})_\tau+\alpha u_\tau=0,\\
&u\cdot\boldsymbol{n}=0.\\
\end{aligned}
\right.
\ee
Here $\al\geq0$ stands for the friction constant which may depend on various elements, such as the property of the boundary and the viscosity of the fluid. When $\al=0$, boundary condition \eqref{NBC} turns to the total Navier-slip boundary \eqref{NSBC}, and when $\al\to\infty$, boundary condition \eqref{NBC} degenerates into the no-slip boundary condition $u\equiv 0$ on the boundary.

The boundary condition \eqref{SLIP} is a special case in a family of boundary conditions proposed by Navier \cite{Navier}, which has been studied extensively in the literature and was attributed to different authors. The boundary condition was called the Navier-Hodge boundary condition in \cite{MM:2009DIE} and the Navier-Lions boundary condition in \cite{PR:2017JFA}. For this reason, we will call it the Navier-Hodge-Lions boundary condition in this paper.

We write $\mD$ to be
 \bes
 \mD=\Sigma\times\bR,
 \ees
 where the cross section $\Sigma\in\bR^2$ is  a unit disc. The domain considered here is a high-degree simplification of the following ``distorted cylinder", i.e.
\[
\t{\mD}=\t{\Sigma}\times\bR,
\]
where $\t{\Sigma}\in\bR^2$ is a simply connected bounded domain with smooth boundary.

Let  $\tilde{\mD}_0$ be a simply connected bounded domain with smooth boundary in $\bR^3$ and $\tilde{\mD}_0\cap \t{\mD}\neq\emptyset$. Existence problem of weak solutions in domain $\t{\mD}_{\text{Union}}:=\t{\mD}\cup\t{\mD}_0$ with the Navier-slip boundary \eqref{NSBC} was addressed in \cite{Konie:2006COLLMATH} and regularity of solutions was also implied there.  On the other hand, if $\t{\mD}_0\subset \t{\mD}$ is an ``obstacle" in $\t{\mD}$, then the two dimensional existence problems and asymptotic behaviors of smooth solutions in domain $\t{\mD}_{\text{Diff}}:=\t{\mD}\backslash \t{\mD}_0$ with the total Navier-slip boundary condition are obtained in \cite{Mucha:2003AAM, Mucha:2003STUDMATH}.

There have also been many pieces of literature in studying the existence, uniqueness and asymptotic behavior of the Navier-Stokes equations in a distorted pipe $\t{\mD}_{\text{Union}}$ or $\t{\mD}_{\text{Diff}}$ with no-slip boundary and with the Poiseuille flow as the asymptotic profile at infinity (Leray's problem: Ladyzhenskaya \cite[p. 77]{Ladyzhenskaya:1959UMN} and \cite[p. 551]{Ladyzhenskaya:1959SPD}). The first remarkable contribution on the solvability of Leray's problem is due to Amick \cite{Amick:1977ASN, Amick:1978NATMA}, who reduced the solvability problem to the resolution of a variational problem related to the stability of the Poiseuille flow in a flat cylinder. However, uniqueness and existence of solutions with large flux are left open. Ladyzhenskaya and Solonnikov \cite{Lady-Sol1980} gave a detailed analysis of this problem on existence, uniqueness and asymptotic behavior of small-flux solutions. One may refer to \cite{AP:1989SIAM, HW:1978SIAM, Pileckas:2002MB} and references for more details on well-posedness, decay and far-field asymptotic analysis of solutions for Leray's problem and related topics. A systematic review and study of Leray's problem can be found in \cite[Chapter XIII]{Galdi:2011SPRINGER}. Recently Wang-Xie in \cite{WX:2019ARXIV} studied uniform structural stability of Poiseuille flows for the 3D axially symmetric solutions in the 3D pipe $\mD$ where a force term  appears on the right hand of equation \eqref{NS}$_1$.

Compared to the no-slip boundary condition, this model with the total Navier-slip or Navier-Hodge-Lions boundary condition has different physical interpretations and gives different mathematical properties. Literature \cite{Mucha:2003AAM, Konie:2006COLLMATH} addressed the existence and regularity problems of weak solutions for the total Navier-slip boundary condition in a pipe, but uniqueness was left open. Also readers can refer to, i.e., \cite{BY:2020ANA, CQ:2010IUMJ, PR:2017JFA, XX:2007CPAM} for some well-posedness results for the Navier-Hodge-Lions boundary condition in different domains.  In this paper, we attempt to derive some uniqueness results for the Navier-Stokes equations with the total Navier-slip or Navier-Hodge-Lions boundary condition in the regular infinite pipe $\mD$. We emphasize that our results below do not require any smallness and decay assumptions.

In this paper, for the boundary condition \eqref{NSBC}, a family of smooth helical solutions will be given, and for the Navier-Hodge-Lions boundary condition \eqref{SLIP}, the trivial shear flow is easy to be discovered. We concern on the characterization and uniqueness of these two types of smooth solutions in $\mD$.

Most of our proof will be carried out in the framework of cylindrical coordinates $(r,\th, z)$ and some of our results are restricted to the axially symmetric solutions. Here we give the formulation of axially symmetric solutions in the cylindrical coordinates which enjoy the following relationship with 3D Euclidian coordinates:
\[
x=(x_1,x_2,x_3)=(r\cos\th,r\sin\th,z).
\]
 A stationary axially symmetric solution of the incompressible Navier-Stokes equations is given as
\[
u=u_r(r,z)\boldsymbol{e_r}+u_{\th}(r,z)\boldsymbol{e_{\th}}+u_z(r,z)\boldsymbol{e_z},
\]
where the basis vectors $\boldsymbol{e_r}$, $\boldsymbol{e_\th}$ ,$\boldsymbol{e_z}$ are
\[
\boldsymbol{e_r}=(\frac{x_1}{r},\frac{x_2}{r},0),\quad \boldsymbol{e_\th}=(-\frac{x_2}{r},\frac{x_1}{r},0),\quad \boldsymbol{e_z}=(0,0,1),
\]
while the components $u_r,u_\th,u_z$, which are independent of $\th$, satisfy
\begin{equation}\label{ASNS}
\left\{
\begin{aligned}
&(u_r\p_r+u_z\p_z)u_r -\frac{(u_\th)^2}{r}+\p_r p=\left(\Delta-\frac{1}{r^2}\right)u_r, \\
&(u_r\p_r+u_z\p_z) u_\th+\frac{u_\th u_r}{r}=\left(\Delta-\frac{1}{r^2}\right)u_\th, \\
&(u_r\p_r+u_z\p_z)u_z+\p_z p=\Delta u_z ,                                    \\
&\na\cdot b=\p_ru_r+\frac{u_r}{r}+\p_zu_z=0,
\end{aligned}
\right.
\end{equation}
where $b=u_r\boldsymbol{e}_r+u_z\boldsymbol{e}_z$.

We can also compute the axi-symmetric vorticity $\o=\nabla\times u=\o_r\boldsymbol{e_r}+\o_\th\boldsymbol{e_\th}+\o_z\boldsymbol{e_z}$  as follows
\[
\o_r=-\p_z u_\th, \ \o_\th=\p_z u_r-\p_r u_z,\  \o_z=\left(\p_r+\frac{1}{r}\right)u_\th,
\]
which satisfies
\be\label{VEQ}
\lt\{
\begin{aligned}
&(u_r\p_r+u_z\p_z)\o_r-\left(\Dl-\f{1}{r^2}\right)\o_r-(\o_r\p_r+\o_z\p_z)u_r=0,\\
&(u_r\p_r+u_z\p_z)\o_\th-\left(\Dl-\f{1}{r^2}\right)\o_\th-\f{u_r}{r}\o_\th-\f{1}{r}\p_z(u_\th)^2=0,\\
&(u_r\p_r+u_z\p_z)\o_z-\Dl \o_z-(\o_r\p_r+\o_z\p_z)u_z=0.
\end{aligned}
\rt.
\ee

In the cylindrical coordinates, the total Navier-slip boundary condition \eqref{NSBC} is represented as
\be\label{NBC1}
\left\{
\begin{aligned}
&\p_ru_\th-\f{u_\th}{r}=0,\\
&\p_ru_z=0,\\
&u_r=0,\\
\end{aligned}
\right.\quad\forall x\in\p \mD,
\ee
while the Navier-Hodge-Lions boundary condition \eqref{SLIP} is given by
\be\label{NBC222}
\left\{
\begin{aligned}
&\p_ru_\th+\f{u_\th}{r}=0,\\
&\p_ru_z=0,\\
&u_r=0,\\
\end{aligned}
\right.\quad\forall x\in\p \mD,
\ee
whose computations are postponed to Appendix \ref{AppA}.

Clearly direct calculation shows that, for arbitrary constants $C_1$ and $C_2$, the  following type of helical solutions
\be\label{helical}
u=C_1r\boldsymbol{e_\th}+C_2 \boldsymbol{e_z}
\ee
solves \eqref{ASNS} with the boundary condition \eqref{NBC1}.
 \begin{figure}[H]
\centering
\includegraphics[scale=0.4]{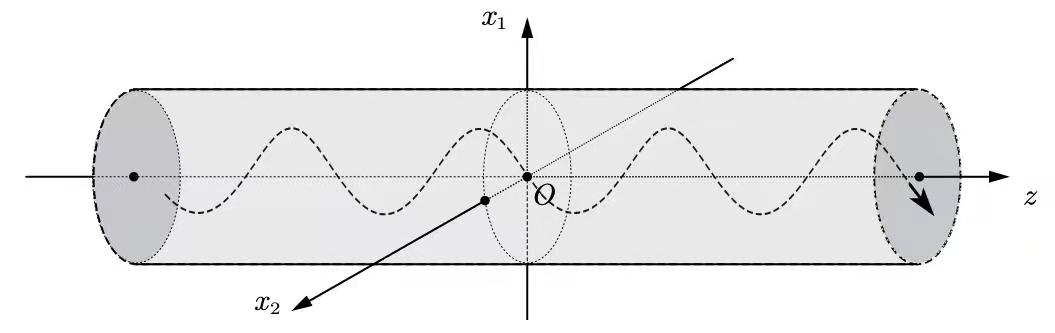}
\caption{A helical solution in the infinite pipe $\mD$}
\end{figure}

 We further note that helical solutions \eqref{helical}, which is smooth in $\mD$, enjoys the following property:

\be\label{bder}\tag{$\ast$}
 \text{\em The solution itself and its gradient are uniformly bounded in } \mD.
\ee

Thus a natural question raises:

\emph{Are helical solutions \eqref{helical} the only smooth solutions of system \eqref{ASNS} with the boundary condition \eqref{NBC1} which enjoys property \eqref{bder}?}

Before answering this question, we recall that the flux $\Phi(z)$ at the cross section $\Sigma$, which is defined by
\[
\Phi(z) := \int_{\Sigma} u(x_h,z)\cdot\boldsymbol{\nu} dx_h,
\]
is a constant. Here $\boldsymbol{\nu}=\boldsymbol{e_z}$ is the unit normal vector of $\Sigma$ pointing to the positive $z$ direction. Actually by using the divergence free condition of the velocity and the boundary condition \eqref{NSBC}$_2$, we have
\bes
\bali
\f{d}{dz}\Phi(z)&=\int_{\Sigma} \f{d}{dz}u_z(x_h,z) dx_h\\
   &=-\int_{\Sigma} (\p_{x_1}u_1+\p_{x_2}u_2)(x_h,z) dx_h \\
   &\xlongequal{Gauss\ formula}-\int_{\p\Sigma} (n_1u_1+n_2u_2)(x_h,z) dS(x_h)\\
   &=-\int_{\p\Sigma}(u\cdot\boldsymbol{n})(x_h,z) dS(x_h)=0,
\eali
\ees
 where $\boldsymbol{n}=(n_1,n_2,0)$ is the unit outer normal vector of $\p\mD$. Then for any  $z\in\mathbb{R}$, we will denote $\Phi(z)=\Phi$.

Our first result in this paper gives a positive answer to the above question in the following two cases:
\begin{itemize}
\item[(i).] the solution is axially symmetric and the flux $\Phi$ is zero (corresponding to $C_2=0$);
\item[(ii).] the solution is $z$-periodic.
\end{itemize}

\begin{theorem}\label{th0}
Let $u$ be a smooth solution of the Navier-Stokes equations in the infinite pipe $\mD$ subject to the total Navier-slip boundary condition \eqref{NSBC}.
 \begin{itemize}
 \item[\bf I.] Suppose $u\in W^{1,\i}(\mD)$ is axially symmetric and $\Phi=0$, then $u$ must be the following type of swirling solutions:
\[
u=C_1r\boldsymbol{e_\th},\q p=\f{C^2_1r^2}{2}.
\]
\item[\bf II.] Suppose $u$ is $z$-periodic, then $u$ must be the following type of helical solutions:
\be\label{helicals}
u=C_1r\boldsymbol{e_\th}+C_2\boldsymbol{e_z},\q p(r,z)={C_1^2r^2}/{2},\q\forall C_1,C_2\in\mathbb{R}.
\ee
\end{itemize}
\end{theorem}

\qed

Besides, we observe that solutions \eqref{helical} enjoy the following property:
\[
\text{\emph{Its swirl component $u_\th$ or vertical component $u_z$ is independent of $z$. }}
\]
In the following theorem, we will conclude that if $u_\th$ or $u_z$  is independent of $z$, then \eqref{helical} are the only group of smooth solutions to \eqref{ASNS} subject to the boundary condition \eqref{NBC1}.

In the case of $u_\th$ being independent of $z$, the bounded $W^{1,\i}$ assumption on the velocity will be extensively relaxed to the following:

\be\label{growass}
\lt\{
\begin{aligned}
&|u_r(r,z)|\leq Cre^{\gamma_0|z|};\\
&|u_z(r,z)|\leq C|z|^{\dl_0};\\
&|\o_\th(r,z)|\leq C|z|^{M_0},
\end{aligned}
\rt.\quad\text{ uniformly with }r\in[0,1],
\ee
for any $\gamma_0<\al\approx3.83171$, $\dl_0<1$, and $M_0>0$. Here $\al$ is the first positive root of the Bessel function $J_1$. We recall that $J_\beta$ are canonical solutions of Bessel's ordinary differential equation
\be\label{BSF}
s^2J_\beta''(s)+sJ_\beta'(s)+(s^2-\beta^2)J_\beta(s)=0,
\ee
which can be expressed by the following series form:
\be\label{BSS}
J_\beta(s)=\sum_{n=0}^\infty\f{(-1)^n}{n!\Gamma(n+\beta+1)}\left(\f{s}{2}\right)^{2n+\beta}.
\ee

In the case of $u_z$ being independent of $z$, no size assumptions such as \eqref{growass} are imposed on the solution.
\begin{remark}
The reason why there is  an $r$ on the righthand of $\eqref{growass}_1$ is that for a smooth solution $u$, in the cylindrical coordinates, $u_r$ vanishes at $r=0$. When doing Taylor expansion of $u_r$ at $r=0$ in the $r$ direction, the zero order derivative term is missing, so it  is reasonable to assume a one order $r$ control on $u_r$ for $r\in[0,1]$.
\end{remark}

\qed

\begin{theorem}\label{th1}
Let $u$ be a smooth solution of the axially symmetric Navier-Stokes equations \eqref{ASNS} in the infinite pipe $\mD$ subject to the total Navier-slip boundary condition \eqref{NSBC}. Then $u$ must be of helical solutions \eqref{helicals} if one of the followings is satisfied.
\begin{itemize}
\item[\bf I.] $u$ satisfies \eqref{growass} and $u_\th$ is independent of $z$-variable;
\item[\bf II.]  $u_z$ is independent of $z$-variable.
\end{itemize}
\end{theorem}

\qed

\begin{remark}
We emphasize that the condition \eqref{growass}$_1$ above is sharp, because we have the following non-trivial counterexample which grows exactly as $Ce^{\al|z|}$ when $z\to\infty$:
\be\label{Count}
\lt\{
\begin{aligned}
&u=-\cosh (\al z)J_1(\al r)\boldsymbol{e_r}+\sinh(\al z)J_0(\al r)\boldsymbol{e_z},\\
&p= -\f{1}{2}\left(\cosh^2(\al z)J_1^2(\al r)+\sinh^2(\al z)J_0^2(\al r)\right).
\end{aligned}
\rt.
\ee
Here $J_0$, $J_1$ are Bessel functions defined in \eqref{BSS}, while $\al\approx3.83171$ is the smallest positive root of $J_1$. One can verify \eqref{Count} being the solution of \eqref{ASNS} with the boundary condition \eqref{NBC1} by direct calculations. Here we leave the details to the interested reader. Unfortunately, our example here can not reflect whether the growing assumptions in \eqref{growass}$_2$ and \eqref{growass}$_3$ are sharp.
\end{remark}

\qed

If we switch the total Navier-slip boundary condition to the Navier-Hodge-Lions boundary condition \eqref{SLIP}, one notices that a non-zero swirling solution does not enjoy \eqref{SLIP}. In this situation, one can derive the following Liouville-type theorem:
\begin{theorem}\label{th2}
Let $u$ be a smooth solution of the Navier-Stokes equations \eqref{NS} in the infinite pipe $\mD$ subject to the Navier-Hodge-Lions boundary condition \eqref{SLIP}. Suppose $\na u$ satisfies
\be\label{COND57}
|\na u(x_h,z)|\leq C|z|^{\beta}
\ee
for some $C>0$ and $0<\beta<1$. Then  $u=\f{\Phi}{\pi}\bl{e_3}$.
\end{theorem}

\qed

There has already been much literature studying Liouville-type results on the Navier-Stokes equations subject to various boundary conditions in various unbounded domains. Readers can refer to \cite{Chae:2014CMP,CW:2016JDE,Seregin:2016NON, Wang:2019JDE,CPZZ:2020ARMA,Pan:2020JMAA} and references therein for more Liouville-type results on the stationary Navier-Stokes equations. Moreover, our results in the above Theorems can be extended from the stationary case to the case of ancient solutions (backward global solutions) under suitable assumptions. However, for simplification of idea presenting, we omit this extension here and leave it to further works. See \cite{GHM:2014CPDE} where the authors established a Liouville-type result for the ancient solution to the Navier-Stokes equations in the half plane with the no-slip boundary condition.

Liouville-type results of ancient solutions is connected to the regularity of solutions to the initial value problem of the non-stationary Navier-Stokes equations. Type I blow-up solutions of the Navier-Stokes initial value problem could not exist provided the Liouville-type result holds for bounded ancient solutions. See \cite{KNSS:2009ACTAMATH,GH:2011CMP}.

Before ending our introduction, we briefly outline our strategy for proofs of Theorem \ref{th0}, Theorem \ref{th1} and Theorem \ref{th2}.  The most important ingredient of proving Theorem \ref{th0} is to show that $\mathbb{S}u \equiv 0$. For {\bf Case I} in Theorem \ref{th0}, we need to first show that $\mathbb{S}u\in L^2(\mD)$, which is not trivial since the pipe is an unbounded domain. In this process, $L^\i$ oscillation boundedness of the pressure in $\mD_{2Z}\backslash \mD_Z$ (see \eqref{ddz} for the definition of $\mD_Z$) is essential, which will be presented in Section \ref{SEC2.2}. Then combining the square integrability of $\mathbb{S}u$ and boundedness of the velocity together with its gradient, a trick of integration by parts and Poincar\'e inequality will indicate that  $u_z$ actually belongs to $L^2(\mD)$, which will result in the vanishing of $\mathbb{S}u$. While vanishing of $\mathbb{S}u$ in {\bf Case II} of Theorem \ref{th0} is directly due to the $x_3$-periodic property by performing the energy estimate in a single period. However, in the process, we must be careful of handing of the pressure, which may not be $x_3$ periodic. After analyzing ingredients of $\mathbb{S}u$, we finally conclude the validity of Theorem \ref{th0}.

The idea for proof of Theorem \ref{th1} is completely different from that of Theorem \ref{th0}. Under the assumption of {\bf Case I} in Theorem \ref{th1}, we will see that the quantity $\O:=\o_\th/r$ satisfies a nice linear elliptic equation with an advection term. Under the growing assumption \eqref{growass} in domain $\mD$, by using the Nash-Moser iteration, we can show that actually $\O\equiv 0$, which indicates that $b=u_r\boldsymbol{e_r}+u_z\boldsymbol{e_z}$ must be harmonic in $\mD$. Then by constructing a barrier function, applying maximum principle and assumptions on $b$, one derives $u_r\equiv0$ and $u_z$ must be a constant. From then on, \eqref{ASNS}{$_2$} is reduced to a linear ordinary differential equation of $u_\th$, and we finally obtain $u_\th=C_1r$. While under the assumption of {\bf Case II} in Theorem \ref{th1}, a rather direct computation implies that $u_r\equiv 0$ by using the divergence-free condition. Then by combining the equations of $u_r$ and $u_z$, independence of $z$ variable for $u_\th$ can be obtained. Then the equation of $u_\th$ will degenerate to an ordinary differential equation with respect to $r$, which will result in $u_\th=C_1r$. At last, trivialness of $u_z$ is achieved by solving a two- dimensional Laplacian equation with Neumann boundary condition.

Proof of Theorem \ref{th2}  is to apply  Lemma \ref{LEM222}, which was originally announced in reference \cite{Lady-Sol1980} as far as the authors know. Denote the energy integral in terms of $v:=u-\f{\Phi}{\pi}\boldsymbol{e}_3$ as follows:
\[
Y(Z):=\int_{Z-1}^Z\int_{\mD_\zeta}|\na v|^2dxd\zeta.
\]
A differential inequality of $Y(Z)$, satisfying the assumption in Lemma \ref{LEM222}, will be derived. In this process, boundary terms coming from integration by parts will be carefully addressed by using the boundary condition, which has a good sign compared with those from the Navier-slip boundary condition.  At last, a direct application of Lemma \ref{LEM222} will imply the vanishing of $Y(Z)$.

For the generalized Navier boundary condition \eqref{NBC} in $\mD$, one can derive that in cylindrical coordinates, \eqref{NBC} is equivalent to
\be\label{gnb}
\left\{
\begin{aligned}
&\p_ru_\th-\f{u_\th}{r}+\al u_\th=0,\\
&\p_ru_z+\al u_z=0,\\
&u_r=0,\\
\end{aligned}
\right.\quad\forall x\in\p \mD.
\ee

For given flux $\Phi:=\int_\Sigma u_z(x_h,z)dx_h=\text{const.}$, we can find a family of bounded smooth solutions satisfying \eqref{ASNS} with boundary condition \eqref{gnb} as follows
\be\label{gnbs}
u=C_1r\chi_{\{\al=0\}}\boldsymbol{e_\th}+\f{2(\al+2)\Phi}{(\al+4)\pi}\lt(1-\f{\al}{\al+2}r^2\rt)\boldsymbol{e_z},\q p={\f{C_1^2r^2}{2}\chi_{\{\al=0\}}}-\f{8\al\Phi}{(\al+4)\pi} z,
\ee
where $C_1$ is an arbitrary constant, and $\chi_{\{\al=0\}}$ is the characteristic function on $\{\al=0\}$, which means
\bes
\chi_{\{\al=0\}}=\lt\{
\bali
&1,\q \al=0,\\
&0, \q \al>0.
\eali
\rt.
\ees

When $\al\rightarrow+\i$, the boundary condition \eqref{gnb} becomes the no-slip boundary and the solution \eqref{gnbs} corresponds to the Hagen-Poisseuille flow in $\mD$. Uniqueness of Hagen-Poisseuille flow is still open for large flux $\Phi$. Our Theorem \ref{th0} states that in the case $\al=0$ and $\Phi=0$, we can show that \eqref{gnbs} are the only bounded smooth solutions of \eqref{ASNS} with the boundary condition \eqref{gnb}. For general $0\leq \al\leq +\i$ and $\Phi$, we have the following conjecture.
\begin{conjecture}\label{CONJ16}
Let $u$ be a smooth solution of the axially symmetric Navier-Stokes equations \eqref{ASNS} in the infinite pipe $\mD$ with the flux $\Phi$ and subject to the Navier-slip boundary condition \eqref{NBC} for any $0\leq \al\leq +\i$. Suppose $u$ and its gradient are uniformly bounded, then the solution $u$ must be of the form \eqref{gnbs}.
\end{conjecture}

\qed

\begin{remark}
From the authors' recent paper \cite{LPY:2022ARXIV}, the uniqueness result there implies that Conjecture \ref{CONJ16} is valid for the case $\al\in (0,+\i)$ if the flux $\Phi$ is small.

\end{remark}

\qed

Throughout this paper, $C_{a,b,c,...}$ denotes a positive constant depending on $a,\,b,\, c,\,...$ which may be different from line to line. For two quantities $A_1$, $A_2$, we denote $A_1\vee A_2=\max\{A_1,\,A_2\}$. Meanwhile, for $Z>1$, we denote
\be\label{ddz}
\mD_\zeta:=\left\{(r,\th,z):\,0\leq r< 1,\,0\leq \th\leq 2\pi,\,-\zeta<z<\zeta\right\},
\ee
the truncated pipe with the length of $2\zeta$. And notations $\mO_\zeta^\pm$ states
\[
\mO_\zeta^+:=\left(\mD_\zeta-\mD_{\zeta-1}\right)\cap\{x\in\mD:\,z>0\},\q\mO_\zeta^-:=\left(\mD_\zeta-\mD_{\zeta-1}\right)\cap\{x\in\mD:\,z<0\},
\]
respectively. We also apply $A\lesssim B$ to denote $A\leq CB$. Moreover, $A\simeq B$ means both $A\lesssim B$ and $B\lesssim A$.

This paper is arranged as follows: Section 2 is devoted to the proof of Theorem \ref{th0}, and the proof of Theorem \ref{th1} could be found in Section 3. Proof of Theorem \ref{th2} will be presented in Section 4.

\section{Proof of Theorem \ref{th0}}\label{sec2}

\q In this section, we devote to proof of Theorem \ref{th0}. Proof of {\bf Case I} is shown in Section \ref{SECI}. In Section \ref{SEC2.1}, we deduce a uniform bound of $\p_z\o_\th$ by using classical energy estimate of \eqref{VEQ}$_2$ and the Moser's iteration. Then it will be applied to derive the $L^\i$ oscillation boundedness of the pressure in Section \ref{SEC2.2}. Based on these preparations, we finish proving {\bf Case I} of Theorem \ref{th0} in Section \ref{SEC2.3}. Proof of {\bf Case II} is directly derived in Section \ref{SECII}.

\subsection{Proof of {\bf Case I}}\label{SECI}

\subsubsection{Uniform bound of ${\p_z\o_\th}$}\label{SEC2.1}

 Denoting $g:=\p_z\o_\th$ and taking $z$-derivative on \eqref{VEQ}$_2$, one arrives
\be\label{EQZVOR}
-\left(\Dl-\f{1}{r^2}\right)g+b\cdot\nabla g=\nabla\cdot F,
\ee
where
\be\label{EPOF}
F:=-\o_\th\p_zb+\left(\f{u_r}{r}\o_\th+2\f{u_\th}{r}\p_zu_\th\right)\boldsymbol{e_z}.
\ee
From \eqref{a3}, we see that $F\in L^\i$ provided $u$ and $\nabla u$ are bounded. Meanwhile, we observe that from the boundary condition \eqref{NBC1}:
\bes
g\equiv 0,\q\text{on}\q\p \mD.
\ees
Now we are ready to state the desired lemma of this section, with its proof based on the Moser's iteration and energy estimate.
\begin{lemma}
Let $(u_r,u_\th,u_z)$ be a smooth solution of \eqref{ASNS} in $\mD$, subject to Navier total slip boundary condition \eqref{NBC1} and $\o_\th$ be the swirl component of its vorticity. Then $\p_z\o_\th$ is uniformly bounded in $\mD$.
\end{lemma}

\pf For $q\geq 1$, we multiply \eqref{EQZVOR} by $qg^{q-1}$ to get
\be\label{TESTT}
-qg^{q-1}\Dl g+\f{q}{r^2}g^q+b\cdot\nabla g^q=qg^{q-1}\nabla\cdot F.
\ee
Noting that
\[
\Dl g^q=\text{div }\left(qg^{q-1}\nabla g\right)=qg^{q-1}\Dl g+q(q-1)g^{q-2}|\nabla g|^2,
\]
one derives from \eqref{TESTT} that
\be\label{SBS}
-\Dl g^q+q(q-1)g^{q-2}|\nabla g|^2+\f{q}{r^2}g^q+b\cdot\nabla g^q= qg^{q-1}\nabla\cdot F.
\ee
Let $\phi$ be a smooth cut-off function in $z$ variable which is bounded up to its second-order derivatives, supported on $[L-1,L+1]$ for some $L\in\mathbb{R}$, which will be specified later. Using $g^q\phi^2$ as a test function to the equation \eqref{SBS} and noting that
\[
q(q-1)\int_{\mathcal{D}}g^{2q-2}|\nabla g|^2\phi^2dx=\f{q-1}{q}\int_{\mathcal{D}}|\nabla g^q|^2\phi^2dx\geq 0,
\]
one deduces
\be\label{ESTG}
\underbrace{\int_{\mD}\nabla g^q\cdot\nabla(g^q\phi^2)dx}_{I_1}+\underbrace{q\int_{\mD}\f{g^{2q}\phi^2}{r^2}dx}_{I_2}+\underbrace{\int_\mD b\cdot\nabla g^q(g^q\phi^2)dx}_{I_3}\leq q\underbrace{\int_\mD g^{2q-1}\nabla\cdot F\phi^2dx}_{I_4}.
\ee
We further denote $f:=g^q$ for convenience. First we see

\[
I_1=\int_{\mD}|\nabla (f\phi)|^2dx-\int_\mD f^2|\nabla\phi|^2dx.
\]
Clearly, $I_2\geq 0$. Using the divergence free property of $b$, one finds $I_3$ satisfies
\[
I_3=\f{1}{2}\int_\mD b\cdot \nabla f^2\phi^2dx= -\int_\mD u_z\p_z\phi^2 f^2dx.
\]
Applying integration by parts, one derives
\[
\begin{split}
I_4=&-q(2q-1)\int_\mD g^{2q-2}\nabla g\cdot F\phi^2dx-q\int_\mD g^{2q-1}F\cdot\nabla\phi^2dx\\
\leq&\f{1}{2}\int_\mD |\nabla(f\phi)|^2dx+Cq^2\int_\mD|F|^2|g|^{2q-2}\phi^2dx+\int_\mD |g|^{2q-1}|F|\phi|\nabla\phi|dx.
\end{split}
\]
Plugging estimates $I_1$--$I_4$ into \eqref{ESTG}, we conclude that
\[
\int_\mD|\nabla(f\phi)|^2dx\leq C\bigg(\|\nabla\phi\|_{L^\infty(\mD)}(\|u_z\|_{L^\infty(\mD)}+\|\nabla\phi\|_{L^\i(\mD)})+q^2\bigg)\int_{\text{supp }\phi}\left(|g|\vee\|F\|_{L^\infty(\mD)}\right)^{2q}dx.
\]
Using the Sobolev imbedding and noting that $\phi$ is supported on an interval whose length equals 2, one arrives
\be\label{ITx}
\begin{split}
\left(\int_{\{x\,:\,\phi=1\}}\left(|g|\vee\|F\|_{L^\i(\mD)}\right)^{6q}dx\right)^{\f{1}{6q}}\leq C^{\f{1}{2q}}\Big(\|\nabla\phi\|_{L^\infty(\mD)}(\|b\|_{L^\infty(\mD)}+\|\nabla\phi\|_{L^\i(\mD)})+q^2\Big)^{\f{1}{2q}}\\
\hskip 5cm\times\left(\int_{\text{supp }\phi}\left(|g|\vee\|F\|_{L^\infty(\mD)}\right)^{2q}dx\right)^{\f{1}{2q}}.
\end{split}
\ee
 Let $\f{1}{2}\leq z_2<z_1\leq 1$ and assume $\phi$ is supported on the interval $[L-z_1,L+z_1]$, and $\phi\equiv 1$ on $[L-z_2,L+z_2]$. Meanwhile, the gradient of $\phi$ satisfies the following estimate:
\[
\|\nabla\phi\|_{L^\infty}\leq\f{C}{z_1-z_2}.
\]
Thus \eqref{ITx} indicates that
\be\label{IT2}
\begin{split}
\left(\int_{\Sigma\times[L-z_2,L+z_2]}\left(|g|\vee\|F\|_{L^\i(\mD)}\right)^{6q}dx\right)^{\f{1}{6q}}&\leq C^{\f{1}{2q}}\Big((z_1-z_2)^{-2}+C_{\|b\|_{L^\infty(\mD)}}+q^2\Big)^{\f{1}{2q}}\\
&\hskip .5cm\times\left(\int_{\Sigma\times[L-z_1,L+z_1]}\left(|g|\vee\|F\|_{L^\infty(\mD)}\right)^{2q}dx\right)^{\f{1}{2q}}.
\end{split}
\ee
Now $\forall k\in\mathbb{N}\cup\{0\}$, we choose $q_k=3^k$ and $z_{1k}=\f{1}{2}+\left(\f{1}{2}\right)^{k+1}$, $z_{2k}=z_{1,k+1}=\f{1}{2}+\left(\f{1}{2}\right)^{k+2}$, respectively. Denoting
\[
\Psi_k:=\left(\int_{\Sigma\times[L-z_{1k},L+z_{1k}]}\left(|g|\vee\|F\|_{L^\infty(\mD)}\right)^{2q_k}dx\right)^{\f{1}{2q_k}},
\]
and iterating \eqref{IT2}, it follows that
\[
\begin{split}
\Psi_{k+1}&\leq C^{\f{1}{2\cdot 3^k}}\left(4^{k+2}+C_{\|b\|_{L^\i}}+3^{2k}\right)^{\f{1}{2\cdot 3^k}}\Psi_{k}\leq\cdot\cdot\cdot\leq \left(C_{\|b\|_{L^\i(\mD)}}\right)^{\f{1}{2}\sum_{j=0}^k3^{-j}}3^{\sum_{j=0}^kj3^{-j}}\Psi_0\leq C_{\|b\|_{L^\i(\mD)}}\Psi_0.
\end{split}
\]
Performing $k\to\i$, the above Moser's iteration implies
\be\label{MOSERR}
\|g\|_{L^\infty(\Sigma\times[L-1/2,L+1/2])}\leq C_{\|b\|_{L^\i(\mD)}}\left(\|g\|_{L^2(\Sigma\times[L-1,L+1])}+\|F\|_{L^\infty(\mathcal{D})}\right).
\ee
Finally, define another cut off function of $z$-variable $\tilde{\phi}$ who has bounded derivatives up to order 2, supported on $[L-2,L+2]$ and $\tilde\phi\equiv1$ in $[L-1,L+1]$. Multiplying \eqref{VEQ}$_2$ by $\o_\th\tilde{\phi}^2$ and integrating on $\mathcal{D}$, one deduces
\[
\int_{\mathcal{D}}|\nabla(\omega_\th\tilde{\phi})|^2dx+\int_{\mathcal{D}}\f{\o_\th^2\tilde{\phi}^2}{r^2}dx=\int_{\mathcal{D}}\o_\th^2|\na\tilde{\phi}|^2dx-\int_{\mathcal{D}}u_z\p_z\tilde{\phi}\o_\th^2\tilde{\phi}dx
-\int_{\mathcal{D}}\f{u_r}{r}\o_\th^2\tilde{\phi}^2dx-2\int_{\mathcal{D}}\f{u_\th}{r}\p_zu_\th\o_\th\tilde{\phi}^2dx.
\]
By the representation of $\nabla u$ \eqref{EDU1}, one derives that
\be\label{EOL2}
\|\nabla \o_\th\|_{L^2(\Sigma\times[L-1,L+1])}\leq C_{\|(u,\na u)\|_{L^\i(\mD)}}.
\ee
 Meanwhile, expression of $F$ \eqref{EPOF} also indicates that
\be\label{EFLI}
\|F\|_{L^\i(\mathcal{D})}\leq C_{\|(u,\na u)\|_{L^\i(\mD)}}.
\ee
Substituting \eqref{EOL2} and \eqref{EFLI} in \eqref{MOSERR}, one concludes that
\[
\|g\|_{L^\infty(\Sigma\times[L-1/2,L+1/2])}\leq C_{\|(u,\nabla u)\|_{L^\i(\mD)}}.
\]
Noting that the right-hand side above is independent of $L$, thus we have derived the uniform boundedness of $g$ in $\mD$.

\qed

\subsubsection{Boundedness of the pressure}\label{SEC2.2}
Based on the boundedness of $\p_z\o_\th$, the $L^\i$ oscillation bound of the pressure $p$ in $D_{2Z}\backslash D_{Z}$ can be obtained. The lemma is stated as follows:
\begin{lemma}\label{lem3.1}
Under the same assumptions of Theorem \ref{th0}, $\forall$ $Z>1$, we have
\be\label{3.111}
\sup_{x\in \mD_{2Z}\backslash\mD_Z}|p(r,z)-p(0,Z)|\leq C,
\ee
where $C>0$ is a uniform constant independent of $Z$.
\end{lemma}
\pf We only consider $\left(\mD_{2Z}\backslash\mD_Z\right)\cap\{x\,:\,z>0\}$ since the rest part is essentially the same. Let us start with the oscillation of the pressure along the $r-$axis. From $\eqref{ASNS}_1$ and the identity
\[
\left(\Dl-\f{1}{r^2}\right)u_r=\p_z\o_\th,
\]
one sees that
\be\label{3.2222}
\p_r p=\p_z\o_\th-(u_r\p_r+u_z\p_z)u_r+\frac{(u_\th)^2}{r}.
\ee
For any $z\in \mathbb{R}$ and $r_1$, $r_2\in[0,1]$, we integrate \eqref{3.2222} with $r$ on $[r_1,r_2]$ to derive
\be\label{EE1}
\begin{split}
p(r_2,z)-p(r_1,z)=&\int_{r_1}^{r_2}\p_z\o_\th dr-\int_{r_1}^{r_2}\left[(u_r\p_r+u_z\p_z)u_r-\f{u^2_\th}{r}\right]dr\\
=&\int_{r_1}^{r_2}\p_z\o_\th(r,z)dr-\f{1}{2}\left(u_r^2(r_2,z)-u_r^2(r_1,z)\right)-\int_{r_1}^{r_2}(u_z\p_zu_r)(r,z)dr\\
&+\int_{r_1}^{r_2}\f{u_\th^2}{r}(r,z)dr.\\
\end{split}
\ee
Noting that
\[
|\nabla u|\simeq|\p_ru_r|+|\p_zu_r|+\left|\f{u_r}{r}\right|+|\p_ru_\th|+|\p_zu_\th|+\left|\f{u_\th}{r}\right|+|\p_ru_z|+|\p_zu_z|,
\]
which follows from \eqref{EDU1}, by the boundedness assumption of $u$ and $\nabla u$, together with the uniform bound of $\p_z\o_\th$ in Section \ref{SEC2.1}, one derives the oscillation bound from \eqref{EE1}:
\be\label{EER}
|p(r_2,z)-p(r_1,z)|\leq\,C(1+\|\p_z\o_\th\|_{L^\infty(\mD_{2Z})})\leq C<\infty,\quad\forall r_1,r_2\in[0,1],\q z\in\mathbb{R},
\ee
where $C$ is an absolute constant which is independent of $r_1$, $r_2$ and $z$. This finishes the oscillation estimate of $p(r,z)$ when $z$ is fixed. Now we turn to the oscillation of $p$ along the $z-$direction. \eqref{ASNS}$_3$ and identity
\[
-\Dl u_z=\f{1}{r}\p_r(r\o_\th)
\]
indicate that
\be\label{3.3333}
\p_zp=-\f{1}{r}\p_r(r\o_\th)-u_r\p_ru_z-u_z\p_zu_z.
\ee
Multiplying \eqref{3.3333} by $r$ and integrating it with respect to $r$ on $(0,1)$, one obtains
\be\label{3.5}
\begin{split}
\f{d}{dz}\int_0^1p(r,z)rdr&=-\un{\int_0^1\p_r(r\o_\th)dr}_{P_1}-\un{\int_0^1(u_r\p_r+u_z\p_z)u_zrdr}_{P_2}.
\end{split}
\ee
Recalling the boundary condition \eqref{NBC1}$_{2,3}$, we find $\o_\th\equiv 0$ on $\p\mD$, which implies $P_1\equiv0$. On the other hand, using the divergence-free condition and integration by parts, we derive
\[
\begin{split}
P_2=&-\int_0^1\p_r\left(ru_r\right)u_zdr+\int_0^1u_z\p_zu_zrdr\\
   =&\int_0^1\p_z\left(ru_z\right)u_zdr+\f{1}{2}\f{d}{dz}\int_0^1u_z^2rdr\\
   =&\f{d}{dz}\int_0^1u_z^2rdr.
\end{split}
\]
\eqref {3.5} indicates
\be\label{3.6}
\f{d}{dz}\int_0^1p(r,z)rdr=-\f{d}{dz}\int_0^1u_z^2(r,z)rdr.
\ee
For any fixed $z\in[Z,2Z]$, we integrate the above indentity from $Z$ to $z$. Then we have
\be\label{EEZ}
\begin{split}
&\left|\int_0^1\big[p(r,z)-p(r,Z)\big]rdr\right|\leq\left|\int_0^1\big[u_z^2(r,z)-u_z^2(r,Z)\big]rdr\right|\leq C.
\end{split}
\ee
Recalling the mean value theorem, there exists $r_*\in[0,1]$ such that
\be\label{EEZ}
|p(r_*,z)-p(r_*,Z)|=\f{\left|\int_0^1\big[p(r,z)-p(r,Z)\big]rdr\right|}{\int_0^1rdr}\leq C.
\ee
This completes the oscillation of $p$ parallel to the $z-$direction. To conclude the general oscillation of the pressure in the pipe, we apply the triangle inequality: for any $r\in[0,1]$, it follows that
\[
|p(r,z)-p(0,Z)|\leq |p(r,z)-p(r_*,z)|+|p(r_*,z)-p(r_*,Z)|+|p(r_*,Z)-p(r,Z)|.
\]
Plugging \eqref{EER} and \eqref{EEZ} into the above inequality, we finally arrive at
\be\label{EEE}
|p(r,z)-p(0,Z)|\leq C,
\ee
where $C$ is an absolute positive constant independent of $r$, $z$ and $Z$. Thus \eqref{3.111} is proved by taking the supremum of \eqref{EEE} over $(r,z)\in[0,1]\times{\big([-2Z,-Z]\cup[Z,2Z]\big)}$.

\qed

\subsubsection{End of the proof}\label{SEC2.3}
In this subsection, we will finish the proof of Theorem \ref{th0}. Namely: \emph{If the flux $\Phi\equiv 0$, any smooth solution of \eqref{ASNS} in an infinite pipe subject to the Navier total slip condition with the velocity and its first-order derivatives being bounded must be a swirling solution}
\[
u=C_1r\boldsymbol{e_\th}.
\]
The proof is divided into three steps: First we show the stress tensor $\mathbb{S}u=\frac{1}{2}\left(\nabla u+(\nabla u)^T\right)$ is globally $L^2$-integrable. Using a 2D Poincar\'e inequality and one insightful observation motivated by \cite{Zhang2021}, we then find that $u_z$ also belongs to $L^2(\mD)$. Finally, we arrive at the vanishing of the stress tensor, which indicates the desired result in Theorem \ref{th0}.

{\bf\noindent $\boldsymbol{L^2}$ boundendness of stress tensor}

Let $\psi:\,\mathbb{R}\to[0,1]$ be a smooth cut-off function satisfying
\[
\psi(l)=\lt\{
\begin{aligned}
&1,\quad l\in[-1,1],\\
&0,\quad |l|\geq2,
\end{aligned}
\rt.
\]
with $\psi'$ and $\psi''$ being bounded. Set
\[
\psi_Z(z):=\psi\left(\f{z}{Z}\right),
\]
where $Z$ is a large positive number. Clearly the derivatives of the scaled cut-off function $\psi_Z$ enjoy
\be\label{EDPSI}
|\p^n_z\psi_Z|\leq \f{C}{Z^n},\q\text{for any }n\in\mathbb{N}.
\ee
Tesing the equation
\[
u\cdot\nabla u+\nabla p=\Dl u
\]
with $u\psi_Z$, we have
\be\label{S2}
\begin{split}
\int_{\mD}\psi_Z u\Dl udx=\int_{\mD}\psi_Zu\Big(u\cdot\nabla u+\nabla\left(p-p(0,Z)\right)\Big)dx.
\end{split}
\ee
To proceed the further calculation in the cylindrical coordinates, we first note that the divergence free property of the velocity indicates
\be\label{4.4}
\begin{split}
\sum_{i,j=1}^3\int_{\mD_{2Z}}\psi_Zu_i\p_{jj}u_idx=\sum_{i,j=1}^3\int_{\mD_{2Z}}\psi_Zu_i\p_{j}\left(\p_ju_i+\p_iu_j\right)dx.\\
\end{split}
\ee
Below, we use the Einstein summation convention for repeated indexes. Using integration by parts, we further derive
\[
\begin{split}
\int_{\mD_{2Z}}\psi_Zu_i\p_{j}\left(\p_ju_i+\p_iu_j\right)dx=&-\underbrace{\int_{\mD_{2Z}}\p_j\psi_{Z}u_i\left(\p_ju_i+\p_iu_j\right)dx}_{T_1}-\underbrace{\int_{\mD_{2Z}}\psi_{Z}\p_ju_i\left(\p_ju_i+\p_iu_j\right)dx}_{T_2}\\
&+\underbrace{\int_{\p \mD_{2Z}}\psi_{Z}u_in_j\left(\p_ju_i+\p_iu_j\right)dS}_{T_3},
\end{split}
\]
where $n_j$ is the $j$-th component of the $\boldsymbol{n}$ -- the unit outward normal vector field on $\p \mD_{2Z}$. Term $T_1$ could be split into two parts, the first half reads
\[
\bali
\int_{\mD_{2Z}}\p_j\psi_{Z}u_i\p_ju_idx=&\f{1}{2}\int_{\mD_{2Z}}\p_j\psi_Z\p_j|u|^2dx=\f{1}{2}\int_{\mD_{2Z}}\p_z\psi_Z\p_z|u|^2dx=-\f{1}{2}\int_{\mD_{2Z}}\p_z^2\psi_{Z}|u|^2dx,
\eali
\]
where we have used the fact that $\psi_{Z}$ is only $z$-dependent and supported in $[-2Z,2Z]$. Similarly, the second half of $T_1$ follows that
\[
\int_{\mD_{2Z}}\p_j\psi_{Z}u_i\p_iu_jdx=\int_{\p \mD_{2Z}}(u\cdot\nabla\psi_Z)(u\cdot\boldsymbol{n})dS-\int_{\mD_{2Z}}\p_z^2\psi_Z u_z^2dx.
\]
Due to the impermeable condition, one sees the first term on the right hand of the above equality is zero. Thus we conclude that
\be\label{TT1}
T_1=-\int_{\mD_{2Z}}\p_z^2\psi_Z\left(\f{1}{2}|u|^2+u_z^2\right)dx.
\ee
Recalling that the stress tensor is defined by
\[
\mathbb{S}u=\f{1}{2}\left(\p_ju_i+\p_iu_j\right)_{1\leq i,j\leq 3},
\]
and using its symmetry, we arrive that
\be\label{TT2}
T_2=\f{1}{2}\sum_{i,j=1}^3\int_{\mD_{2Z}}\psi_Z\left(\p_ju_i+\p_iu_j\right)^2dx=2\int_{\mD_{2Z}}\psi_Z|\mathbb{S}u|^2dx.
\ee
Now applying the Navier-slip condition \eqref{NSBC}$_1$, one notes that
\[
n_j\left(\p_ju_i+\p_iu_j\right)=c(x)n_i,
\]
where $c(x)$ is a scalar-valued function. Inserting this identity to $T_3$, we find
\be\label{TT3}
T_3=\int_{\p \mD_{2Z}}c\psi_Z\left(u\cdot\boldsymbol{n}\right)dS=0.
\ee
Next we come back to the right hand side of \eqref{S2}. Noting $u$ is divergence-free, integration by parts shows
\be\label{TT4}
\begin{split}
&\int_{\mD_{2Z}}u\psi_Z\Big(u\cdot\nabla u+\nabla\big[p-p(0,Z)\big]\Big)dx=\int_{\mD_{2Z}}\psi_Z u_i\p_i\left(\f{1}{2}|u|^2+\big[p-p(0,Z)\big]\right)dx\\
=&\underbrace{\int_{\p \mD_{2Z}}\psi_Z(u\cdot\boldsymbol{n})\Big(\f{1}{2}|u|^2+\big[p-p(0,Z)\big]\Big)dS}_{T_4}-\int_{\mD_{2Z}}\p_z\psi_Z u_z\left(\f{1}{2}|u|^2+\big[p-p(0,Z)\big]\right)dx.
\end{split}
\ee
Here $T_4$ above also vanishes by the stationary wall condition $\eqref{NBC1}_3$. Therefore we arrive that by plugging \eqref{TT1}, \eqref{TT2}, \eqref{TT3}, \eqref{TT4} into \eqref{S2}
\be\label{ELAST}
2\int_{\mD_{2Z}}\psi_Z|\mathbb{S}u|^2dx=\int_{\mD_{2Z}}\p^2_{z}\psi_Z\left(\f{1}{2}|u|^2+u_z^2\right)dx+\underbrace{\int_{\mD_{2Z}}\p_z\psi_Z u_z\left(\f{1}{2}|u|^2+\big[p-p(0,Z)\big]\right)dx}_{T_5}.
\ee
Recalling \eqref{EDPSI}, the bounds on the derivatives of scaled cut-off function $\psi_Z$, and the boundedness of $u$ and pressure, one derives from \eqref{ELAST} that
\[
\int_{\mD_{2Z}}\psi_Z|\mathbb{S}u|^2dx\leq C|\mD_{2Z}|\left(Z^{-2}+Z^{-1}\right)\leq C,
\]
where $C$ is a universal constant depending only on the $L^\i$ bound of $u$ and $\nabla u$ given in the assumption. After letting $Z\to\infty$, the above inequality shows the stress tensor is globally $L^2$-integrable:
\be\label{ESTRESS}
\int_{\mD}|\mathbb{S}u|^2dx\leq C<\infty.
\ee

{\bf \noindent$\boldsymbol{L^2}$ boundedness of $u_z$}\label{sec2.3.2}

First we observe that $\|u_z\|_{L^2(\mD)}$ can be controlled by $\|\p_ru_z\|_{L^2(\mD)}$ under the assumption that the flux $\Phi\equiv 0$.  Noting that
\bes
\f{1}{|\Sigma|}\int_\Sigma u_z(x_h,z)dx_h=\f{1}{|\Sigma|}\Phi=0,
\ees
then we apply the one dimensional Poincar\'e inequality to derive
\[
\bali
\int_{\Sigma}|u_z(r,z)|^2dx_h=&\int_\Sigma\left|u_z(x_h,z)-\f{1}{|\Sigma|}\int_\Sigma u_z(x_h,z)dx_h\right|^2dx_h\\
                            \leq& S^2_0\int_\Sigma|\na_h u_z(x_h,z)|^2dx_h=S^2_0\int_\Sigma|\p_r u_z(r,z)|^2dx_h,
\eali
\]
where $\na_h=(\p_1,\p_2)$ and $S_0$ is independent of $z\in\mathbb{R}$. Integrating with $z$-variable on $\mathbb{R}$, we arrive
\be\label{EUZ}
\|u_z\|_{L^2(\mD)}\leq S_0\|\p_r u_z\|_{L^2(\mD)}.
\ee
However, we cannot get the $L^2$ boundedness of $\p_ru_z$ directly from \eqref{ESTRESS}. In fact, by the expression of the stress tensor \eqref{STEN}, one only has the uniform $L^2$ bound of $(\p_zu_r+\p_ru_z)$. Nevertheless, one observes
\[
\begin{split}
\int_{\mD_{2Z}}\left(\p_ru_z\right)^2dx&=\int_{\mD_{2Z}}\left(\p_zu_r+\p_ru_z\right)^2dx-\int_{\mD_{2Z}}\left(\p_zu_r\right)^2dx-2\int_{\mD_{2Z}}\p_ru_z\p_zu_rdx\\
&\leq C+2\Big|\underbrace{\int_{\mD_{2Z}}\p_ru_z\p_zu_rdx}_{T_6}\Big|.
\end{split}
\]
Now it remains to derive the boundedness of $T_6$. With idea motivated by \cite{Zhang2021}, after using the divergence free of $u$ and integration by parts, we deduce
\begin{small}
\[
\begin{split}
&\q\int_{D_{2Z}}\p_ru_z\p_zu_rdx\\
&=-2\pi\int_{-2Z}^{2Z}\int_0^1u_z\p^2_{rz}(ru_r)drdz=2\pi\int_{-2Z}^{2Z}\int_0^1u_z\p^2_z(ru_z)drdz\\
&=-\underbrace{\int_{D_{2Z}}\left(\p_zu_z\right)^2dx}_{T_7}+2\pi\underbrace{\left(\int_{0}^1u_z(r,2Z)\p_zu_z(r,2Z)rdr-\int_{0}^1u_z(r,-2Z)\p_zu_z(r,-2Z)rdr\right)}_{T_8}.
\end{split}
\]
\end{small}
Here $T_7$ can be bounded by the $L^2$ norm of stress tensor \eqref{ESTRESS}, while $T_8$ is controlled by the $L^\i$ bounds of $u$ and $\nabla u$. Noting that $T_6$ is estimated uniformly with respect to $Z$. This, together with \eqref{EUZ} implies
\bes
\|u_z\|_{L^2(\mD)}\leq C<\infty.
\ees
{\bf\noindent Vanishing of $\boldsymbol{\int_{\mD}|\mathbb{S}u|^2}$ and finishing of the proof}

Based on the $L^2$ bound of $u_z$, now we can estimate $T_5$ in \eqref{ELAST} in an alternative approach, by using H\"older inequality:
\[
|T_5|\leq\sup_{x\in \mD_{2Z}\backslash\mD_{Z}}\left|\f{1}{2}|u|^2+\big[p-p(0,Z)\big]\right|\f{C}{Z}\|u_z\|_{L^2(\mD_{2Z})}|\mD_{2Z}|^{1/2}\leq CZ^{-1/2}.
\]
Thus we deduce from \eqref{ELAST}
\[
\int_{\mD_{2Z}}\psi_Z|\mathbb{S}u|^2dx\leq C|\mD_{2Z}|Z^{-2}+CZ^{-1/2}\rightarrow 0,\q \text{as } Z\rightarrow +\i,
\]
which indicates that
\be\label{zero}
\int_{\mD}|\mathbb{S}u|^2dx=0
\ee
by letting $Z\to\infty$. By the expression of $\mathbb{S}u$ \eqref{STEN}, one finds
\[
u_r\equiv\p_z u_\th\equiv\p_zu_z\equiv\p_ru_z\equiv0,\q\p_ru_\th=\f{u_\th}{r}.
\]
The above estimates, together with the vanishing flux ($\Phi=0$), indicate
\[
u_z\equiv0, \q\text{and}\q u_\th=Cr.
\]
Thus we conclude that $u=Cr\boldsymbol{e_\th}$, which is a swirling solution.

\qed

\subsection{Proof of Case II }\label{SECII}

If $u$ is $z$-periodic with the minimal period $L>0$, then we can drop the restriction $\Phi=0$ in {\bf Case I}.\\
The proof is straightforward. Set $w=u-\f{\Phi}{\pi}\boldsymbol{e_z}$, then $(w,p)$ satisfies the following system:
\begin{equation}\label{NS-r}
\begin{cases}
(w+\f{\Phi}{\pi}\boldsymbol{e_z})\cdot\na w+\na p-\Dl w=0\,,\quad&\text{in $\mD$}\,,\\
\na\cdot w=0\,,&\text{in $\mD$}\,,\\
w(x_h,z)=w(x_h,z+L)\,,&\text{in $\mD$}\,,\\
(\mathbb{S}w\cdot\boldsymbol{n})_\tau=0\,,\quad
w\cdot\boldsymbol{n}=0\,,&\text{on $\p \mD$}\,.
\end{cases}
\end{equation}

We first {\bf claim} that, the pressure $p(x_h,z)$ has the following decomposition:
\begin{equation}\label{p-dec}
	p(x_h,z) =\,az +\tilde p (x_h,z),
\end{equation}
where $a$ is a constant, and  $\tilde p$  is $z$-periodic with the minimal period $L>0$. Set
\begin{equation*}
	a_0(x_h)=\f{1}{L}\int_0^L(\p_zp)(x_h,z)dz\,,
\end{equation*}
we decompose $\p_zp$ as
\begin{equation}\label{pzp}
	(\p_zp)(x_h,z)=a_0(x_h)+\left((\p_zp)(x_h,z)-a_0(x_h)\right):=a_0(x_h)+p_1(x_h,z)\,.
\end{equation}
By using the equation and $z-$periodicity of the solution, it is easy to check that $p_1(x_h,z)$ is $L$-periodic with respect to $z$ and that
\begin{equation}\label{p10}
	\int_0^Lp_1(x_h,\tilde{z})\,d\tilde{z}=0\,.
\end{equation}
Integrating \eqref{pzp} with $z$ on $[0,z]$, one derives that
\begin{equation*}
p(x_h,z)=p(x_h,0)+a_0(x_h)z+\int_0^zp_1(x_h,\tilde{z})\,d\tilde{z}\,.
\end{equation*}
It is worth noting that $\int_0^zp_1(x_h,\tilde{z})\,d\tilde{z}$ is periodic in the $z$-direction due to \eqref{p10} and periodicity of $p_1$. Hence,
\begin{equation*}
\tilde{p}(x_h,z)=p(x_h,0)+\int_0^zp_1(x_h,\tilde{z})\,d\tilde{z}
\end{equation*}
is periodic in the $z$-direction, and
\begin{equation*}
p(x_h,z)=a_0(x_h)z+\tilde{p}(x_h,z)\,.
\end{equation*}
Finally, also from the equations and $z$-periodicity of the solution, we deduce $\nabla_{h}p=(\nabla_{h}a_0)z+(\nabla_{h}\tilde{p})$ is periodic with respect to $z$. Thus we get
\[
a_0(x_h)=\text{constant}:=a
\]
since $\nabla_ha_0$ must be zero. Therefore we conclude
\begin{equation*}
p(x_h,z)=az+\tilde{p}(x_h,z)\,.
\end{equation*}
This finishes the proof of the {\bf claim}. 

Next, multiplying $w$ on both sides of \eqref{NS-r}$_1$, and integrating on $\Sigma\times[0,L]$, one has
\[\int_{\Sigma\times[0,L]}w\cdot\Delta w\,dx=\int_{\Sigma\times[0,L]}w\cdot\left(\big(w+\f{\Phi}{\pi}\boldsymbol{e_z}\big)\cdot\na w+\na p\right)\,dx\,.\]
It follows from \eqref{NS-r}$_2$--\eqref{NS-r}$_4$ that
\[\int_{\Sigma\times[0,L]}w\cdot\Delta w\,dx=-2\int_{\Sigma\times[0,L]}|\mathbb{S}w|^2\,dx,\]
where we have used the technique from \eqref{4.4} to \eqref{TT2} and $z$-periodicity of $w$ to do integration by parts. There are no boundary terms generated. Also 
\[\int_{\Sigma\times[0,L]}w\cdot\left(\big(w+\f{\Phi}{\pi}\boldsymbol{e_z}\big)\cdot\na w+\na p\right)\,dx=0\,,\]
where we have used the decomposition \eqref{p-dec} and $\int_{\Sigma}w_3\,dx_h=0$ to deal with the pressure term.
Hence, we have that
\[\int_{\Sigma\times[0,L]}|\mathbb{S}w|^2\,dx=0\,,\]
which deduces that $\mathbb{S}w=0$. It is well known that if $\mathbb{S}w=0$, then $w$ has the form $w=Ax+B$ (see \cite[\S 6]{KO:1988RMS}), where $A$ is a skew-symmetric matrix with constant entries and $B$ is a constant vector, that is,
\begin{equation*}
	w=
	\begin{pmatrix}
		0 & -a_1 &-a_2 \\
		a_1 & 0 & -a_3 \\
		a_2 & a_3 & 0
	\end{pmatrix}
	\begin{pmatrix}
		x_1 \\
		x_2 \\
		z
	\end{pmatrix}
	+\begin{pmatrix}
		b_1 \\
		b_2 \\
		b_3
	\end{pmatrix}=
	\begin{pmatrix}
		-a_1\,x_2-a_2\,z+b_1	 \\
		a_1\,x_1-a_3\,z+b_2	\\
		a_2x_1+a_3x_2+b_3
	\end{pmatrix}\,,
\end{equation*}
where $a_i\,,b_i$ ($i=1\,,2\,,3$) are some constants. Note that $w$ is periodic with respect to $z$, we have that $a_2=a_3=0$. Next, from $\int_{\Sigma}w_3\,dx_h=0$, we obtain that $b_3=0$. Finally, note that $w_r=0$ on $r=1$, we get that $b_1=b_2=0$. Thus, we have proved that
\begin{equation*}
	w=
	\begin{pmatrix}
		-a_1\,x_2	 \\
		a_1\,x_1\\
		0
	\end{pmatrix}
	=a_1r\,\boldsymbol{e_\theta}\,.
\end{equation*}
Therefore, $u=w+\f{\Phi}{\pi}\boldsymbol{e_z}=a_1r\,\boldsymbol{e_\theta}+\f{\Phi}{\pi}\boldsymbol{e_z}$.

\qed

Let us give some discussions of Theorem \ref{th0} here. Based on our previous proof in this section, we naturally believe that if the vanishing of $\Phi$ is abandoned, then an axially symmetric solution must be a helical solution:
\be\label{HELI}
u=C_1r\boldsymbol{e_\th}+C_2\boldsymbol{e_z}
\ee
 even without the $z$-periodic condition. However, our method in this paper fails when we handle solutions with the flux $\Phi\neq 0$, because we can no longer apply the Poincar\'e inequality in Section \ref{sec2.3.2} to derive the $L^2$ integrability of $u_z$. Meanwhile, if we denote
\[
c_0:=\f{1}{|\Sigma|}\int_\Sigma u_z(x_h,z)dx_h=\f{1}{|\Sigma|}\Phi,
\]
then $u_z-c_0$ enjoys a similar Poincar\'e inequality as \eqref{EUZ}:
\[
\|u_z-c_0\|_{L^2(\mD)}\leq S_0\|\p_ru_z\|_{L^2(\mD)},
\]
which guarantees the $L^2$ boundedness of $u_z-c_0$. However, one additional term appears in $T_5$ of \eqref{ELAST}, which is:
\[
T_5':=c_0\int_{\mD_{2Z}}\p_z\psi_Z\left(\f{1}{2}|u|^2+\big[p-p(0,Z)\big]\right)dx.
\]
Without any integrability of the head pressure $\f{1}{2}|u|^2+\big[p-p(0,Z)\big]$, we can only show $T_5'$ is bounded, which results in
\bes
\int_\mD|\mathbb{S}u|^2dx<C<\i.
\ees
However, we are unable to conclude $T_5'\to 0$ as $Z\to\infty$, thus vanishing of $\int_\mD|\mathbb{S}u|^2dx$ can not be obtained. In fact, using integration by parts on $z$ in $T_5'$, we have
\bes
T_5'=-c_0\int_{\mD_{2Z}}\psi_Z\p_z\left(\f{1}{2}|u|^2+p\right)dx.
\ees
By following the argument in Section \ref{sec2}, one derives
\[
\int_\mD|\mathbb{S}u|^2dx=-\lim_{Z\to\infty}\f{c_0}{2}\int_\mD\psi_Z\p_z\left(\f{1}{2}|u|^2+p\right)dx
\]
instead of \eqref{zero}. Recalling \eqref{3.6}, one deduces that
\be\label{Pre}
\int_\mD|\mathbb{S}u|^2dx=-\lim_{Z\to\infty}\f{c_0}{4}\int_\mD\psi_Z\p_z\left(u_r^2+u_\th^2-u_z^2\right)dx.
\ee
Thus if $\p_z\left(u_r^2+u_\th^2-u_z^2\right)\in L^1(\mD)$ (or $\p_z\left(u_r^2+u_\th^2-u_z^2\right)$ has a fixed sign), one concludes the following identity by Lebesgue's dominated convergence theorem (or monotone convergence theorem):
\be\label{GEQ}
\int_\mD|\mathbb{S}u|^2dx+\f{c_0}{4}\int_\mD\p_z\left(u_r^2+u_\th^2-u_z^2\right)dx=0.
\ee

At the moment, even with identities \eqref{Pre} and \eqref{GEQ} for bounded (up to first-order derivatives) smooth axisymmetric solutions of stationary Navier-Stokes equations in $\mD$ subject to the total Navier-slip boundary condition in hand, we neither show the trivialness of $\mathbb{S}u$, nor find a nontrivial bounded solution apart from \eqref{HELI} which satisfies conditions of Theorem \ref{th0}. Indeed, we leave characterization of the non-zero flux solutions  in Conjecture \ref{CONJ16}.

Nevertheless, a direct observation of \eqref{GEQ} indicates that: If $u$ is independent of $z$, then the right hand side of \eqref{GEQ} vanishes and we can conclude $\mathbb{S}u\equiv 0$, and thus conclude that $u=C_1r\boldsymbol{e_\th}+C_2\boldsymbol{e_z}$ as we desire. In the next section, we see that only $u_\th$ or $u_z$ being independent of $z$ is adequate for us to derive Theorem \ref{th1}. Besides, the asymptotic assumptions of $u$ and its derivatives can be largely loosened.

\section{Proof of Theorem \ref{th1}}

\subsection{Proof of Case I}
\q Let us outline the proof at the beginning of this section: Under the assumptions of {\bf Case I} in Theorem \ref{th1}, our first step is showing $\o_\th\equiv0$, which indicates $b=u_r\boldsymbol{e_r}+u_z\boldsymbol{e_z}$ must be harmonic in $\mD$. Then by applying the boundary  condition and the asymptotic behavior of $b$, one derives $u_r\equiv0$ and $u_z$ must be a constant. From then on \eqref{ASNS}$_2$ turns to a linear ordinary differential equation of $u_\th$, and we finally prove $u_\th=C_1r$.

\subsubsection{Vanishing of ${\o_\th}$}
Noting that $u_\th$ is independent of $z$, we find \eqref{VEQ}$_2$ now turns to
\[
(u_r\p_r+u_z\p_z)\o_\th-\left(\Dl-\f{1}{r^2}\right)\o_\th-\f{u_r}{r}\o_\th=0.
\]
From the Navier-slip boundary condition \eqref{NBC1}, one has
\[
\o_\th=\p_zu_r-\p_ru_z=0,\quad\text{on }\p \mD.
\]
Denoting $\O:=\f{\o_\th}{r}$, direct calculation shows
\be\label{IBO}
\left\{
\begin{split}
&(u_r\p_r+u_z\p_z)\O-\left(\Dl+\f{2}{r}\p_r\right)\O=0,\quad&\text{in }\quad \mD;\\
&\O=0,\quad&\text{on }\q\p \mD.
\end{split}
\right.
\ee
In the following, we first provide a mean value inequality of $\O$ deduced by Moser's iteration.
\begin{lemma}
Assume $b=u_r\boldsymbol{e_r}+u_z\boldsymbol{e_z}$ is a smooth divergence-free axially symmetric vector field. Then any weak solution $\O$ of boundary value problem \eqref{IBO} satisfies the following mean value inequality:
\be\label{MEQ}
\sup_{x\in \mD_{\tau_2Z}}|\O|\leq C_q(\tau_1-\tau_2)^{-\f{q}{q-2}}\left(1+\|u_z\|_{L^\infty(\mD_{Z}\backslash \mD_{Z/2})}\right)^{\f{q}{q-2}}Z^{-\f{q}{q-2}}\left(\int_{\mD_{\tau_1Z}}|\O|^{2}dx\right)^{\f{1}{2}},
\ee
for any $q>2$, $Z>1$, and $\f{1}{2}\leq\tau_2<\tau_1\leq1$.
\end{lemma}

\pf We only prove \eqref{MEQ} with $\tau_1=1$, $\tau_2=\f{1}{2}$ for simplicity, since the general case could be derived by a direct scaling strategy. For any real number $l\geq1$, we find $h:=\O^l$ satisfies
\be\label{Omega}
\Dl h-l(l-1)\O^{l-2}|\na\O|^2+\f{2}{r}\p_r h-b\cdot\nabla h=0.
\ee
Set $\f{1}{2}\leq\sigma_2<\sigma_1\leq1$ and choose $\zeta=\zeta(z)$ to be a smooth cut-off function satisfying
\[
\left\{\begin{array}{l}
\operatorname{supp} \zeta \subset [-\sigma_1,\sigma_1], \quad \zeta=1 \quad \text { in } [-\sigma_2,\sigma_2],\\
 0 \leqslant \zeta \leqslant 1,\\
\left|\zeta'\right| \lesssim\frac{1}{\sigma_{1}-\sigma_{2}}.
\end{array}\right.
\]
Denoting $\zeta_Z(z):=\zeta\left(\f{z}{Z}\right)$ and testing \eqref{Omega} with $\zeta_Z^2 h$, noting that
\[
l(l-1)\int_{\mD_{\sigma_1Z}}\O^{2l-2}|\nabla \O|^2\zeta_Z^2dx=\f{l-1}{l}\int_{\mD_{\sigma_1Z}}|\nabla \O^l|^2\zeta_Z^2dx\geq 0,
\]
we arrive
\be\label{EM0}
\un{\int_{\mD_{\sigma_1Z}}\Dl h\zeta_Z^2hdx}_{M_1}+\un{\int_{\mD_{\sigma_1Z}}\f{2}{r}\p_rh\zeta_Z^2hdx}_{M_2}-\un{\int_{\mD_{\sigma_1Z}}b\cdot\nabla h\zeta_Z^2hdx}_{M_3}\geq0.
\ee
Next we handle $M_1$--$M_3$ term by term. Using integration by parts and direct calculations, we first see
\be\label{EM1}
\begin{split}
M_1=-\int_{\mD_{\sigma_1Z}}\nabla h\cdot\nabla(\zeta_Z^2h)dx=-\int_{\mD_{\sigma_1Z}}|\na(h\zeta_Z)|^2dx+\int_{\mD_{\sigma_1Z}}h^2|\zeta_Z'|^2dx.
\end{split}
\ee
Here the boundary term of the cylindrical surface is cancelled because $h=0$ on $\p \mD$, while those coming from the cross sections $D\cap\{z=\pm \sigma_1Z\}$ vanish due to the cut off function $\zeta_Z$ is compactly supported. On the other hand, using axisymmetry of the solution
\be\label{EM2}
M_2=2\pi\int_{\bR}\int_0^1\p_r(h^2\zeta_Z^2)drdz=-2\pi\int_{\bR}h^2(0,z)\zeta^2_Z(z)dz\leq0.
\ee
Before we bound $M_3$, let us introduce the \emph{stream function} of axisymmetric velocity field $b=u_r\boldsymbol{e_r}+u_z\boldsymbol{e_z}$. By the divergence-free property $\p_r(ru_r)+\p_z(ru_z)=0$, there exists a scalar function $L_\theta=L_\theta(r,z)$ such that
\be\label{psii3}
-\p_zL_\theta=u_r,\quad\text{and}\quad\frac{1}{r}\p_r(rL_\theta)=u_z.
\ee
Using integration by parts again together with boundary condition $h=0$ on $\p \mD$, we derive that
\[
\begin{split}
M_3&=\f{1}{2}\int_{D_{\sigma_1Z}}b\cdot\nabla h^2\zeta_Z^2dx=-\int_{D_{\sigma_1Z}}u_z\zeta_Z\zeta_Z'h^2dx=-2\pi\int_{\bR}\int_0^1\p_r(rL_\th)\zeta_Z\zeta_Z'h^2drdz\\
&=4\pi\int_{\bR}\int_0^1(rL_\th)\p_r(h\zeta_Z)h\zeta_Z'drdz.
\end{split}
\]
By the mean value theorem and \eqref{psii3}, there exists $\tilde{r}\in(0,r)$ such that
\[
rL_\th(r,z)=\tilde{r}u_z(\tilde{r},z)r,
\]
thus we can further bound $M_3$ by
\be\label{EM3}
\begin{split}
|M_3|&\leq4\pi\|u_z\|_{L^\infty(\mD_{\sigma_1Z}\backslash\mD_{\sigma_2Z})}\int_{\bR}\int_0^1|\nabla (h\zeta_Z)h\zeta_Z'|rdrdz\\
&\leq\f{1}{2}\int_{\mD_{\sigma_1Z}}|\na(h\zeta_Z)|^2dx+2\|u_z\|_{L^\infty(\mD_{\sigma_1Z}\backslash \mD_{\sigma_2Z})}^2\int_{\mD_{\sigma_1Z}}h^2|\zeta_Z'|^2dx.
\end{split}
\ee
Now substituting \eqref{EM1}, \eqref{EM2}, and \eqref{EM3} in \eqref{EM0}, taking the maximum of $\zeta_Z'$, it follows that
\be\label{EEE1}
\int_{\mD_{\sigma_1Z}}|\na(h\zeta_Z)|^2dx+2\pi\int_{\bR}h^2(0,z)\zeta^2_Z(z)dz\leq\f{C\left(1+\|u_z\|_{L^\infty(\mD_{\sigma_1Z}\backslash \mD_{\sigma_2Z})}^2\right)}{(\si_1-\si_2)^2Z^2}\int_{\mD_{\sigma_1Z}}h^2dx.
\ee
Recalling $h=0$ on $\p \mD$, for any fixed $z\in\mathbb{R}$, the following 2D Poincar\'e inequality holds:
\[
\|h(\cdot,z)\zeta_Z(z)\|^2_{L^2(\Sigma)}\leq C\left\|\p_r\left[h(\cdot,z)\zeta_Z(z)\right]\right\|^2_{L^2(\Sigma)},
\]
where $C>0$ here is independent of $z$. Integrating with $z$ on $\bR$ and taking the square root, one has the following 3D Poincar\'e inequality
\be\label{Po}
\|h\zeta_Z\|_{L^2(\mD_{\sigma_1Z})}\leq C\|\p_r(h\zeta_Z)\|_{L^2(\mD_{\sigma_1Z})}.
\ee
For any $q\in(2,6)$, Interpolation, Sobolev inequality and \eqref{Po} imply that
\be\label{EEE2}
\begin{split}
\|h\zeta_Z\|_{L^q(\mD_{\sigma_1Z})}&\leq \|h\zeta_Z\|^{s}_{L^6(\mD_{\sigma_1Z})}\|h\zeta_Z\|^{1-s}_{L^2(\mD_{\sigma_1Z})}\leq C\|\nabla(h\zeta_Z)\|^{s}_{L^2(\mD_{\sigma_1Z})}\|h\zeta_Z\|^{1-s}_{L^2(\mD_{\sigma_1Z})}\\
&\leq C\|\nabla(h\zeta_Z)\|^{s}_{L^2(\mD_{\sigma_1Z})}\|\p_r(h\zeta_Z)\|^{1-s}_{L^2(\mD_{\sigma_1Z})}\leq C \|\nabla(h\zeta_Z)\|_{L^2(\mD_{\sigma_1Z})}.
\end{split}
\ee
Here $s\in(0,1)$ depends on $q$. Combining \eqref{EEE1} and \eqref{EEE2}, we derive
\[
\|h\|_{L^q(\mD_{\sigma_2Z})}\leq\f{C\left(1+\|u_z\|_{L^\infty(\mD_{\sigma_1Z}\backslash \mD_{\sigma_2Z})}\right)}{(\sigma_1-\sigma_2)Z}\|h\|_{L^2(\mD_{\sigma_1Z})},
\]
which is equivalent to
\be\label{PIT}
\left(\int_{D_{\sigma_2Z}}|\O|^{ql}dx\right)^{\f{1}{ql}}\leq\f{C^{1/l}\left(1+\|u_z\|_{L^\infty(D_{\sigma_1Z}\backslash D_{\sigma_2Z})}\right)^{1/l}}{(\si_1-\si_2)^{\f{1}{l}}Z^{\f{1}{l}}}\left(\int_{D_{\sigma_1Z}}|\O|^{2l}dx\right)^{\f{1}{2l}}.
\ee
Now for any $k=0,1,2,...$, we choose $l_k=\left(\f{q}{2}\right)^k$ and $\si_{1k}=\f{1}{2}+\left(\f{1}{2}\right)^{k+1}$, $\si_{2k}=\f{1}{2}+\left(\f{1}{2}\right)^{k+2}$, respectively. Denoting
\[
\Psi_k:=\left(\int_{D_{\sigma_{1k}Z}}|\O|^{2l_k}dx\right)^{\f{1}{2l_k}},
\]
and noting that
\[
D_{\sigma_{1k}Z}\backslash D_{\sigma_{2k}Z}\subset D_Z\backslash D_{Z/2},\q\forall k=0,1,2,...,
\]
then \eqref{PIT} follows that
\be\label{IT1}
\begin{split}
\Psi_{k+1}&\leq C^{\left(\f{2}{q}\right)^k}2^{k\left(\f{2}{q}\right)^k}\left(1+\|u_z\|_{L^\infty(\mD_{Z}\backslash \mD_{Z/2})}\right)^{\left(\f{2}{q}\right)^k}Z^{-\left(\f{2}{q}\right)^{k}}\Psi_{k}\\[4mm]
&\leq\cdot\cdot\cdot\\[4mm]
&\leq C^{\sum_{j=0}^k\left(\f{2}{q}\right)^j}2^{\sum_{j=0}^k j\left(\f{2}{q}\right)^j}\left(1+\|u_z\|_{L^\infty(\mD_{Z}\backslash \mD_{Z/2})}\right)^{\sum_{j=0}^k\left(\f{2}{q}\right)^j}Z^{-\sum_{j=0}^k\left(\f{2}{q}\right)^{j}}\Psi_0.
\end{split}
\ee
Performing $k\to\infty$, then iteration \eqref{IT1} implies a mean value inequality of $\O$, that is
\[
\sup_{x\in \mD_{Z/2}}|\O|\leq C_q\left(1+\|u_z\|_{L^\infty(\mD_{Z}\backslash \mD_{Z/2})}\right)^{\f{q}{q-2}}Z^{-\f{q}{q-2}}\left(\int_{\mD_{Z}}|\O|^{2}dx\right)^{\f{1}{2}},
\]
for any $q>2$. This completes the proof of Lemma \ref{lem3.1}.

\qed

Since $u_z$  satisfies \eqref{growass}$_2$ in $\mD_{Z}$, \eqref{MEQ} indicates that
\be\label{IT2}
\sup_{x\in \mD_{\tau_2 Z}}|\O|^2\leq C_q(\tau_1-\tau_2)^{-\f{2q}{q-2}} Z^{\f{2(\dl_0-1)q}{q-2}}\int_{\mD_{\tau_1 Z}}|\O|^{2}dx.
\ee
However, due to the lack of boundedness of the second-order derivatives of $u$, we are unable to control $\|\O\|_{L^2(\mD_Z)}$ at the moment. Next we will use an algebraic trick to convert the $L^2$-norm on the right hand side of \eqref{IT2} to  an $L^1$-norm. This trick comes from Li-Schoen \cite{Li-Schoen1984}. Here goes the lemma:
\begin{lemma}[modified mean value inequality]\label{lem3.2}
Suppose $b=u_r\boldsymbol{e_r}+u_z\boldsymbol{e_z}$ is a smooth divergence-free axisymmetric vector field and $\|u_z\|_{L^\infty(D_Z)}\les Z^{\dl_0}.$ Then any weak solution $\O$ of boundary value problem \eqref{IBO} satisfies the following mean value inequality for any $q>2$, $Z>1$:
\be\label{IT4}
\sup_{x\in \mD_{Z/2}}|\O|\leq C_q Z^{\f{(\dl_0-1)q}{q-2}}\int_{\mD_Z}|\O|dx.
\ee
\end{lemma}
\pf For any $\f{1}{2}\leq\tau_2<\tau_1\leq1$, \eqref{IT2} implies that
\[
\sup_{x\in \mD_{\tau_2 Z}}|\O|^2\leq C_q(\tau_1-\tau_2)^{-\f{2q}{q-2}} Z^{\f{2(\dl_0-1)q}{q-2}} \left(\sup_{x\in \mD_{\tau_1 Z}}|\O|^2\right)^{1/2}\int_{D_{Z}}|\O|dx.
\]
Denoting $\tau_{1k}=\tau_{2,k+1}=1-\left(\f{1}{2}\right)^{k+2}$, $\tau_{2k}=1-\left(\f{1}{2}\right)^{k+1}$, and $\Phi_k:=\sup_{x\in D_{\tau_{2k}Z}}|\O|^2$, it follows that
\be\label{IT3}
\Phi_k\leq C_q2^{\f{2qk}{q-2}} Z^{\f{2(\dl_0-1)q}{q-2}}\Phi_{k+1}^{1/2}\int_{D_{Z}}|\O|dx.
\ee
Iterating \eqref{IT3} from $k=0$ to infinity, one arrives
\[
\begin{split}
\sup_{x\in \mD_{Z/2}}|\O|^2&\leq C_q^{\sum_{j=0}^\infty 2^{-j}}2^{\f{2q}{q-2}\sum_{j={0}}^\infty\f{j}{2^j}}\left( Z^{\f{2(\dl_0-1)q}{q-2}}\right)^{\sum_{j=0}^\infty 2^{-j}}\left(\int_{\mD_{Z}}|\O|dx\right)^{\sum_{j=0}^\infty 2^{-j}}\\
&\leq C_qZ^{\f{2(\dl_0-1)q}{q-2}}\left(\int_{\mD_{Z}}|\O|dx\right)^2,
\end{split}
\]
which follows that
\[
\sup_{x\in \mD_{Z/2}}|\O|\leq C_q Z^{\f{(\dl_0-1)q}{q-2}}\int_{\mD_Z}|\O|dx.
\]
This completes the proof of Lemma \ref{lem3.2}.

\qed

Finally, one notes that
\[
\int_{\mD_Z}|\O|dx\leq 2\pi\|\o_\th\|_{L^\i(\mD_{Z})}\int_{-Z}^Z\int_0^1drdz\les  Z^{M_0+1}.
\]
Therefore, as long as $\o_\th$ is of polynomial growth (see \eqref{growass}$_3$) when $z\to\infty$, we can infer from \eqref{IT4} that
\be\label{IT5}
\sup_{x\in \mD_{Z/2}}|\O|\leq C_q Z^{\f{(\dl_0-1)q}{q-2}+1+M_0}.
\ee

For any fixed $\dl_0<1$ and $M_0>0$, we can always choose $q$ which is larger than but close enough to $2$ such that \eqref{IT5} indicates
\[
\sup_{x\in \mD_{Z/2}}|\O|\les Z^{-\gamma}
\]
for some $\gamma>0$. This proves $\o_\th$ vanishes in $\mD$ by performing $Z\to\infty$.

\subsubsection{Vanishing of $u_r$ and constancy of $u_z$}\label{SEC3.2}
Noting that $\nabla\times b=\o_\th\boldsymbol{e_\th}\equiv0$ and the divergence-free property of $b$, we apply the Lagrange's formula for del to deduce
\[
-\Dl b=\nabla\times\nabla\times b-\nabla(\text{div }b)=0,
\]
which indicates
\[
\left(\Dl-\f{1}{r^2}\right)u_r=0;\q\q\Dl u_z=0.
\]
To prove vanishing of $u_r$, for $\dl>0$ being small, we consider the auxiliary function $\eta_\dl$ which is defined by
\[
\eta_\dl(x):=J_1\big((\al-\dl)r\big)\cosh\big((\al-\dl)z\big).
\]
Here $J_1$ is the Bessel function which is defined in \eqref{BSS} and satisfies \eqref{BSF} with $\beta=1$, while $\al$ is the smallest positive root of $J_1$. Direct calculation shows
\[
\left(\Dl-\f{1}{r^2}\right)\eta_\delta=\left(\p^2_r+\f{1}{r}\p_r+\p^2_z-\f{1}{r^2}\right)\eta_\dl=0.
\]
Owing to $u_r$ is growing as \eqref{growass}$_1$, we choose $\dl<<1$ small enough such that $\g_0<\al-2\dl$. Using the concavity of $J_1((\al-\dl)r)$  on the subset of $\{r:\,0\leq r\leq 1\}$ where $J_1((\al-\dl)r)$ is increasing, one has
\[
\eta_\dl\geq J_1(\al-\dl)r\cosh\big((\al-\dl)z\big)\geq C_\dl re^{(\gamma_0+\dl)|z|},
\]
where $C_\dl>0$ is a constant depends only on $\dl$. Then the condition \eqref{growass}$_1$ indicates that
\[
\lim_{|z|\to\infty}\f{|u_r(r,z)|}{\eta_\delta(r,z)}=0, \text{ uniformly with }r=\sqrt{x_1^2+x_2^2}\in[0,1].
\]
Therefore, for any fixed $\ve>0$ and $\dl$, there exists an $N_{\ve,\delta}>0$ such that
\[
\lt\{
\begin{aligned}
&\left(\Delta-\f{1}{r^2}\right)(\ve\eta_\delta\pm u_r)=0,\quad\forall x\in \mD_M,\\
&\ve\eta_\delta\pm u_r\geq0,\quad\forall x\in\partial \mD_M =\big[\p\mD\cap\{-M\leq z\leq M\}\big]\cup\big[\mD\cap\{z=\pm M\}\big],
\end{aligned}
\rt.
\]
for any $M>N_{\ve,\delta}$. The maximum principle indicates
\be\label{5.9}
|u_r(x)|\leq \ve\eta_\delta(x),\quad\forall x\in \mD_M.
\ee
By performing $M\to\infty$, one finds the estimate \eqref{5.9} actually holds for all $x\in \mD$. Thus $u_r\equiv0$ is proved by the arbitrariness of $\ve>0$.

Finally, the divergence-free of $u$ implies $\p_zu_z=-\f{1}{r}\p_r(ru_r)\equiv 0$ in $\mD$. The vanishing of $\o_\th$ and $u_r$ indicates $\p_ru_z\equiv0$. Thus $u_z$ must be a constant. This consequently indicates
\be\label{Solb}
b=C_2\boldsymbol{e_z}
\ee
for some constant $C_2\in\bR$.

\subsubsection{End of the proof}\label{SEC313}
Now substituting \eqref{Solb} in \eqref{ASNS}$_2$ and noting that $u_\th$ is independent of $r$, one arrives the following ODE of $u_\th$
\[
u_\th''(r)+\f{1}{r}u_\th'(r)-\f{1}{r^2}u_\th(r)=0.
\]
This ODE, which is of Eulerian type, is solved by
\[
u_\th(r)=\f{C_0}{r}+C_1r,
\]
for any $C_0$, $C_1\in\bR$. Smoothness of $u_\th$ forces that $C_0=0$. Thus we conclude that
\[
u=u_\th\boldsymbol{e_\th}+b=C_1r\boldsymbol{e_\th}+C_2\boldsymbol{e_z},
\]
which completes the proof of Theorem \ref{th1}.

\qed

\begin{remark}
Unlike Theorem \ref{th0}, Theorem \ref{th1} actually needs weaker assumptions,\eqref{growass}, on the boundedness of solutions.
%
As stated in the introduction, assumption \eqref{growass}$_1$ is sharp due to the non-trivial counterexamples in \eqref{Count} which grow no slower than $Ce^{\al |z|}$ as $z\to\infty$.
%
Meanwhile, the counterexample in \eqref{Count} has zero vorticity and zero flux in the cross section $\Sigma$. Identities \eqref{Pre} and \eqref{GEQ} no longer hold for the solution in \eqref{Count} since we have no boundedness of the head pressure $H:=\f{1}{2}|u|^2+p-p(0,Z)$ in $\mD_{2Z}\backslash \mD_Z$.

\end{remark}

\qed

\subsection{Proof of Case II}
Actually, if $u_z$ is independent of $z$-variable, by the divergence free condition \eqref{ASNS}$_{4}$, we have 
\[
\p_r(ru_r)=-r\p_z u_z=0,
\]
which indicates that $ru_r=f(z)$ for some smooth function $f(z)$. Then using the boundary condition \eqref{NBC1}$_{3}$, we deduce that $f(z)=0$, which implies 
\[u_r=0.\]
Next, it follows from \eqref{ASNS}$_{1}$ and \eqref{ASNS}$_{3}$ that
\[\p_z(u_{\theta})^2=r\p_{r}\p_{z}p=r\p_r\left(\p^2_{r}u_z+{\f{1}{r}\p_r u_z}\right):=g(r)\,,\]
which deduces that
\[(u_{\theta})^2(r,z)=g(r)z+(u_{\theta})^2(r,0)\,.\]
 Note that $g(r)\equiv 0$ for any $r\in[0,1]$. Actually, if for some $r_0$ such that $g(r_0)\neq0$, we can obtain a contradiction by taking
 \[
 z=\frac{-1-(u_{\theta})^2(r_0,0)}{g(r_0)}.
 \]
 Thus, we have that $u_{\theta}$ is independent of $z$-variable. Following the argument in Section \ref{SEC313}, one concludes that $u_\th=C_1 r$.

 Then if we go back to the \eqref{ASNS}$_1$, we see that $\p_r p=C^2_1r$, which indicates that
 \[
 p=\f{1}{2} C^2_1r^2+f(z),
 \]
 for some smooth function $f(z)$. From \eqref{ASNS}$_3$, by using that $u_z$ is independent of $z$, we can have that $\p_zp=f'(z)=C$ for some constant $C$. At last we see that $u_z$ satisfies the following two dimensional Laplacian equation in $\Sigma$ with Neumann boundary condition:
 \be\label{neu}
 \lt\{
 \bali
 &\Dl_h u_z=C,\\
 &\p_{\boldsymbol{n}}u_z=0.
 \eali
 \rt.
 \ee
Integrating directly on $\Sigma$ for \eqref{neu}$_1$ and using the boundary condition, we can obtain $C=0$. Then multiplying  \eqref{neu}$_1$ by $u_z$ and integrating on $\Sigma$, we can have $\na_h u_z=0$, which implies that $u_z \equiv C_2$ for some constant $C_2$.

\qed

\section{Proof of Theorem \ref{th2}}
In this section we derive the proof of Theorem \ref{th2}, which shows a solution to the Navier-Stokes equations \eqref{NS} with the $Navier-Hodge-Lions$ boundary condition \eqref{SLIP} in the pipe $\mD$ must be a parallel flow $\f{\Phi}{\pi}\bl{e_z}$, without axisymmetric assumptions. Our method is motivated by \cite{Lady-Sol1980} in which authors show the uniqueness result for problems with the homogeneous Dirichlet boundary condition. Before proving the theorem, we need the following lemma that states the asymptotic behavior of a function satisfies an ordinary differential inequality.
\begin{lemma}\label{LEM222}
Let $Y(\zeta)\nequiv 0$ be a nondecreasing nonnegative differentiable function satisfying
\[
Y(\zeta)\leq\Psi(Y'(\zeta)),\q\forall\zeta>0.
\]
Here $\Psi:\,[0,\infty)\to[0,\infty)$ is a monotonically increasing function with $\Psi(0)=0$ and there exists $C,\,\tau_1>0$, $m>1$, such that
\[
\Psi(\tau)\leq C\tau^m,\q\forall\tau>\tau_1.
\]
Then
\[
\liminf_{\zeta\to+\infty}\zeta^{-\f{m}{m-1}}Y(\zeta)>0.
\]
\end{lemma}

\qed

The next lemma on solving the divergence problem in a truncated pipe will be applied to bound the term related to pressure in the further proof:
\begin{lemma}[See \cite{Bme1, Bme2}, also \cite{Galdi:2011SPRINGER}, Chapter III]\label{LEMTrun}
Let $D=\Sigma\times [0,1]$, $f\in L^2(D)$ with
\[
\int_D fdx=0,
\]
then there exists a vector valued function ${V}:\,D\to\mathbb{R}^3$ belongs to $H^1_0(D)$ such that
\be\label{LEM2.11}
\nabla\cdot {V}=f,\q\text{and}\q\|\nabla {V}\|_{L^2(D)}\leq C\|f\|_{L^2(D)}.
\ee
Here $C>0$ is an absolute constant.
\end{lemma}

\qed

The following lemma gives a Poincar\'e inequality for vectors in $H^1(\Sigma)$ with only vanishing normal direction on the boundary. Readers can find some hint of the proof in Galdi \cite[Page 71, Exercise II.5.6]{Galdi:2011SPRINGER}.
\begin{lemma}\label{P057}
Let ${f}=f_1\boldsymbol{e_1}+f_2\boldsymbol{e_2}$ be a two dimensional vector function with components in $H^1(\Sigma)$, and ${f}\cdot\bar{\boldsymbol{n}}=0$ on $\partial \Sigma$, where $\bar{\boldsymbol{n}}$ is the unit outer normal of $\p \Sigma$. Then the following Poincar\'e inequality holds
\[
\|{f}\|_{L^2(\Sigma)}\leq C\|\nabla_h {f}\|_{L^2(\Sigma)},
\]
where $\na_h=(\p_{x_1},\p_{x_2})$ is the gradient operator on $x_1$ and $x_2$ direction.
\end{lemma}

\qed

\leftline{\textbf{Proof of Theorem \ref{th2}:} Denoting $v:=u-\f{\Phi}{\pi}\bl{e_z}$, one deduces from \eqref{NS} that}
\be\label{57subt}
v\cdot\na v+\f{\Phi}{\pi}\p_zv+\na p-\Dl v=0.
\ee
We multiply \eqref{57subt} by $v$ and integrate on $\mD_\zeta$, it follows that
\be\label{57Mair}
\un{\int_{\mathcal{D}_\zeta}v\cdot\Dl vdx}_{LHS}=\int_{\mathcal{D}_\zeta}v\left(v\cdot\nabla v+\f{\Phi}{\pi}\p_zv+\nabla p\right)dx.
\ee
Using integration by parts, the left hand side of \eqref{57Mair} follows that
\[
\begin{split}
LHS=-\int_{\mD_\zeta}|\na v|^2dx+\un{\f{1}{2}\int_{\p\mD_\zeta}\f{\p|v|^2}{\p\bl{n}}dS}_{B_1},
\end{split}
\]
where $\bl{n}$ is the unit outer normal vector on $\p\mD_\zeta$. In the cylindrical coordinate, one writes
\[
v=v_r\bl{e_r}+v_\th\bl{e_\th}+v_z\bl{e_z}.
\]
Thus one has
\be\label{B1157}
\begin{split}
B_1=&\un{\int_{\p\mD_\zeta\cap\p\mD}\big(v_r\p_rv_r+v_\th\p_rv_\th+v_z\p_rv_z\big)dS}_{B_{11}}+{\int_{\mD\cap\{z=\zeta\}}\big(v_r\p_zv_r+v_\th\p_zv_\th+v_z\p_zv_z\big)dx_h}\\
&-{\int_{\mD\cap\{z=-\zeta\}}\big(v_r\p_zv_r+v_\th\p_zv_\th+v_z\p_zv_z\big)dx_h}.
\end{split}
\ee
Since $v$ satisfies the boundary condition \eqref{NBC222}, one derives
\[
B_{11}=-\int_{\p\mD_\zeta\cap\p\mD}v_\th^2dS.
\]
Substituting above in \eqref{B1157}, one derives
\be\label{ELHS}
LHS\leq -\int_{\mD_\zeta}|\na v|^2dx-\int_{\p\mD_\zeta\cap\p\mD}v_\th^2dS+\int_{\mD\cap\{z=\pm\zeta\}}|v||\p_zv|dx_h.
\ee
Now we turn to the right hand side of \eqref{57Mair}. Integrating by parts, one deduces that
\be\label{Maint357}
\begin{split}
\int_{\mathcal{D}_\zeta}v\big(v\cdot\nabla v+\nabla p\big)dx=&\int_{\mD\cap\{z=\zeta\}}v_3\left(\f{1}{2}|v|^2+p\right)dx_h-\int_{\mD\cap\{z=-\zeta\}}v_3\left(\f{1}{2}|v|^2+p\right)dx_h.
\end{split}
\ee
Meanwhile, one notices
\be\label{Maint457}
\int_{\mD_\zeta}\f{\Phi}{\pi} v\cdot\p_zvdx=\f{\Phi}{2\pi}\left(\int_{\mD\cap\{z=\zeta\}}|v|^2dx_h-\int_{\mD\cap\{z=-\zeta\}}|v|^2dx_h\right).
\ee
Substituting \eqref{ELHS}, \eqref{Maint357} and \eqref{Maint457} in \eqref{57Mair}, one arrives at
\[
\begin{split}
\int_{\mD_\zeta}|\nabla v|^2dx\leq C\int_{\mD\cap\{z=\pm\zeta\}}|v|(|\na v|+|v|+|v|^2)dx_h-\int_{\mD\cap\{z=\zeta\}}v_3pdx_h+\int_{\mD\cap\{z=-\zeta\}}v_3pdx_h.
\end{split}
\]
Integrating with $\zeta$ on $[Z-1,Z]$, where $Z\geq 1$, it follows that
\be\label{ET057}
\int_{Z-1}^Z\int_{\mD_\zeta}|\nabla v|^2dxd\zeta\leq C\Bigg(\un{\int_{\mO^{+}_Z\cup\mO^{-}_Z}|v|(|\na v|+|v|+|v|^2)dx}_{T_1}+\un{\Big|\int_{\mO^{+}_Z\cup\mO^{-}_Z}v_3pdx\Big|}_{T_2}\Bigg).
\ee
In the following we only work on integrations on $\mO^{+}_Z$ since the rest part are similar. Using Cauchy-Schwarz inequality and Gagliardo-Nirenberg inequality, one deduces
\[
T_1\leq C\left(\|v\|_{L^2(\mO^+_Z)}\|\na v\|_{L^2(\mO^+_Z)}+\|v\|_{L^2(\mO^+_Z)}^2+\|v\|^{3/2}_{L^2(\mO^+_Z)}\|\na v\|^{3/2}_{L^2(\mO^+_Z)}+\|v\|_{L^2(\mO_Z^+)}^3\right).
\]
Applying Lemma \ref{P057} in each cross section of the pipe, one finds
\be\label{ET157}
T_1\leq C\left(\|\na v\|^2_{L^2(\mO^+_Z)}+\|\na v\|^3_{L^2(\mO^+_Z)}\right).
\ee
Now it remains to bound the pressure term $T_2$ in \eqref{ET057}. Noticing that
\[
\int_{\mD\cap\{x_3=z\}}v_3(x_h,z)dx_h\equiv0,\q\forall z\in\mathbb{R},
\]
we deduce that
\[
\int_{\mO_Z^+}v_3dx=0,\q\forall Z\geq 1.
\]
Using Lemma \ref{LEMTrun}, one derives the existence of a vector field $V$ satisfying \eqref{LEM2.11} with $f=v_3$. By the momentum equation \eqref{NS}$_1$, one arrives
\[
\int_{\mO_Z^+}v_3 pdx=-\int_{\mO_Z^+}\nabla p\cdot {V}dx=\int_{\mO_Z^+}\left(v\cdot\na v+\f{\Phi}{\pi}\p_zv-\Dl v\right)\cdot {V}dx.
\]
integration by parts, one derives
\[
\int_{\mO_Z^+}v_3 pdx=\sum_{i,j=1}^3\int_{\mO_Z^+}\left(\p_iv_j-v_iv_j-\f{\Phi}{\pi}\delta_{i3}v_j\right)\p_iV_jdx.
\]
Here $\delta_{ij}$ is the Kronecker symbol. By applying H\"older inequality and \eqref{LEM2.11} in Lemma \ref{LEMTrun}, one deduces that
\[
\left|\int_{\mO_Z^+}v_3 pdx\right|\leq C\left(\|\na v\|_{L^2(\mO^+_Z)}+\|v\|_{L^4(\mO^+_Z)}^2+\Phi\|v\|_{L^2(\mO^+_Z)}\right)\|v_3\|_{L^2(\mO^+_Z)}.
\]
Similarly as we bound $T_1$, using Lemma \ref{P057} and Gagliardo-Nirenberg inequality, one concludes
\be\label{ET257}
T_2\leq C\left(\|\nabla v\|_{L^2(\mO^+_Z)}^2+\|\nabla v\|_{L^2(\mO^+_Z)}^3\right).
\ee
Substituting \eqref{ET157} and \eqref{ET257} (together with their related estimates on $\mO^-_Z$ ), one deduces
\[
\int_{Z-1}^Z\int_{\mD_\zeta}|\na v|^2dxd\zeta\leq C\left(\|\nabla v\|_{L^2(\mO^+_Z\cup \mO^-_Z)}^2+\|\nabla v\|_{L^2(\mO^+_Z\cup\mO^-_Z)}^3\right),\q\forall Z\geq 1.
\]
Therefore, letting
\[
Y(Z):=\int_{Z-1}^Z\int_{\mD_\zeta}|\na v|^2dxd\zeta,
\]
one deduces
\[
Y(Z)\leq C\left(Y'(Z)+\left(Y'(Z)\right)^{3/2}\right).
\]
Applying Lemma \ref{LEM222}, one concludes
\[
Y(Z)=\int_{Z-1}^Z\int_{\mD_\zeta}|\na v|^2dxd\zeta\geq C_0Z^3,\q\forall Z\geq 1
\]
for some $C_0>0$. However, this creates a paradox since $|\na u|$ satisfies \eqref{COND57}. This proves $u=\f{\Phi}{\pi}\bl{e_3}$.

\qed

\begin{appendix}
\titleformat{\section}[block]{\Large\center\sc}{\sc Appendix}{0.5em}{}[]

\section{Computation of the boundary condition}\label{AppA}

\q Here we give a derivation of the boundary condition \eqref{NSBC} and \eqref{SLIP} in the cylindrical coordinates. First, noting that
\be\label{A1}
0=u\cdot\boldsymbol{n}=u_r.
\ee

In cylindrical coordinates, the gradient operator is represented by
\[
\na= \boldsymbol{e_r}\p_r +\boldsymbol{e_\th}\f{\p_\th}{r}+\boldsymbol{e_z}\p_z.
\]
Then we can calculate the matrix $\na u$ in cylindrical coordinates and write it as a form of tensor product as follows
\be\label{EDU1}
\begin{split}
\nabla u&=\p_ru_r\boldsymbol{e_r}\otimes\boldsymbol{e_r}+{\lt(\f{1}{r}\p_\th u_r-\frac{u_\th}{r}\rt)}\boldsymbol{e_r}\otimes\boldsymbol{e_\th}+\p_zu_r\boldsymbol{e_r}\otimes\boldsymbol{e_z}\\
&+\p_ru_\th\boldsymbol{e_\th}\otimes\boldsymbol{e_r}+{\lt(\f{1}{r}\p_\th u_\th+\frac{u_r}{r}\rt)}\boldsymbol{e_\th}\otimes\boldsymbol{e_\th}++\p_zu_\th\boldsymbol{e_\th}\otimes\boldsymbol{e_z}\\
&+\p_ru_z\boldsymbol{e_z}\otimes\boldsymbol{e_r}+\f{1}{r}\p_\th u_z\boldsymbol{e_z}\otimes\boldsymbol{e_\th}+\p_zu_z\boldsymbol{e_z}\otimes\boldsymbol{e_z}.
\end{split}
\ee
Equivalently
\be\label{a3}
\bali
\na  u=&
\begin{pmatrix}
\p_r u_r & \f{1}{r}\p_\th u_r-\frac{1}{r}u_\th &\p_z u_r \\
\p_r u_\th &\f{1}{r}\p_\th u_\th+\frac{1}{r}u_r & \p_z u_\th \\
\p_r u_z & \f{1}{r}\p_\th u_z & \p_z u_z
\end{pmatrix}
:
\begin{pmatrix}
\boldsymbol{e_r}\otimes\boldsymbol{e_r} & \boldsymbol{e_r}\otimes\boldsymbol{e_\th} &\boldsymbol{e_r}\otimes\boldsymbol{e_z} \\
\boldsymbol{e_\th}\otimes\boldsymbol{e_r}  &\boldsymbol{e_\th}\otimes\boldsymbol{e_\th} & \boldsymbol{e_\th}\otimes\boldsymbol{e_z} \\
\boldsymbol{e_z}\otimes\boldsymbol{e_r} & \boldsymbol{e_z}\otimes\boldsymbol{e_\th} & \boldsymbol{e_z}\otimes\boldsymbol{e_z}
\end{pmatrix}\\
  :=&\begin{pmatrix}
\p_r u_r & \f{1}{r}\p_\th u_r-\frac{1}{r}u_\th &\p_z u_r \\
\p_r u_\th &\f{1}{r}\p_\th u_\th+\frac{1}{r}u_r & \p_z u_\th \\
\p_r u_z & \f{1}{r}\p_\th u_z & \p_z u_z
\end{pmatrix}
:\mathcal{A}.
\eali
\ee
Also a direct computation shows that
\[
\na\times u=\left(\f{1}{r}\p_\th u_z-\p_zu_\th\right)\bl{e_r}+\left(\p_zu_r-\p_ru_z\right)\bl{e_\th}+\f{1}{r}\left(\p_r(ru_\th)-\p_\th u_r\right)\bl{e_z}.
\]
Then ${\mathbb  S} u$ under the base $\mathcal{A}$ is represented by
\be\label{STEN}
{\mathbb  S} u=\begin{pmatrix}
\p_r u_r & \f{1}{2}\lt(\f{1}{r}\p_\th u_r+\p_r u_\th-\f{1}{r}u_\th\rt) &\f{1}{2}(\p_z u_r+\p_r u_z) \\
 \f{1}{2}\lt(\f{1}{r}\p_\th u_r+\p_r u_\th-\f{1}{r}u_\th\rt) &\f{1}{r}\p_\th u_\th+\frac{1}{r}u_r & \f{1}{2}\lt(\p_z u_\th+\f{1}{r}\p_\th u_z\rt)\\
\f{1}{2}(\p_z u_r+\p_r u_z) & \f{1}{2}\lt(\p_z u_\th+\f{1}{r}\p_\th u_z\rt) & \p_z u_z
\end{pmatrix}:\mathcal{A}.
\ee
 Since the outward normal vector $\boldsymbol{n}=\boldsymbol{e_r}$, we have
\[
{\mathbb  S} u\cdot \boldsymbol{n}=\p_r u_r\boldsymbol{e_r}+\f{1}{2}\lt(\f{1}{r}\p_\th u_r+\p_r u_\th-\f{1}{r}u_\th\rt)\boldsymbol{e_\th}+\f{1}{2}(\p_z u_r+\p_r u_z)\boldsymbol{e_z}.
\]
Then in cylinder coordinates, one has
\[
({\mathbb  S} u\cdot \boldsymbol{n})_{\tau}= \f{1}{2}\lt(\f{1}{r}\p_\th u_r+\p_r u_\th-\f{1}{r}u_\th\rt)\boldsymbol{e_\th}+\f{1}{2}(\p_z u_r+\p_r u_z)\boldsymbol{e_z},
\]
and
\[
\na\times u\times \boldsymbol{n}=-\left(\p_zu_r-\p_ru_z\right)\bl{e_z}+\f{1}{r}\left(\p_r(ru_\th)-\p_\th u_r\right)\bl{e_\th}.
\]
This, together with \eqref{A1}, the boundary condition \eqref{NSBC} and \eqref{SLIP} in the cylindrical coordinates read

\[
\left\{
\begin{aligned}
&\p_ru_\th-\f{u_\th}{r}=0,\\
&\p_ru_z=0,\\
&u_r=0,\\
\end{aligned}
\right.\quad\forall x\in\p \mD,
\]
and
\[
\left\{
\begin{aligned}
&\p_ru_\th+\f{u_\th}{r}=0,\\
&\p_ru_z=0,\\
&u_r=0,\\
\end{aligned}
\right.\quad\forall x\in\p \mD,
\]
where we have used the fact that $\p_\th u_r=\p_z u_r=0$ on the boundary.
\end{appendix}

\section*{Acknowledgments}
\addcontentsline{toc}{section}{Acknowledgments}
\q\ The authors wish to thank Professors Qi S. Zhang, Xin Yang and Na Zhao, and Mr. Chulan Zeng for helpful discussions and proof reading. Z. Li is supported by Natural Science Foundation of Jiangsu Province (No. BK20200803) and National Natural Science Foundation of China (No. 12001285). X. Pan is supported by National Natural Science Foundation of China (No. 11801268, 12031006). J. Yang is supported by National Natural Science Foundation of China (No. 12001429).
\section*{Conflict of interest statement}

On behalf of all authors, the corresponding author states that there is no conflict of interest.

\section*{Data availability statement}

Data sharing not applicable to this article as no datasets were generated or analysed during the current study.

\medskip

 {\footnotesize

{\sc Z. Li: School of Mathematics and Statistics, Nanjing University of Information Science and Technology, Nanjing 210044, China}

  {\it E-mail address:}  zijinli@nuist.edu.cn

\medskip

 {\sc X. Pan: College of Mathematics, Nanjing University of Aeronautics and Astronautics, Nanjing 211106, China}

  {\it E-mail address:}  xinghong\_87@nuaa.edu.cn

\medskip
{\sc J. Yang: School of Mathematics and Statistics, Northwestern Polytechnical University, Xi'an 710129, China}

  {\it E-mail address:} yjqmath@nwpu.edu.cn
}

\begin{thebibliography}{}
\addcontentsline{toc}{section}{References}
\small
\setlength{\itemsep}{-3 pt}

\bibitem{Amick:1977ASN} {\sc C. J. Amick}: Steady solutions of the Navier-Stokes equations in unbounded channels and pipes. {\it Ann. Scuola Norm. Sup. Pisa Cl. Sci.} (4) 4 (1977), no. 3, 473--513.

\bibitem{Amick:1978NATMA} {\sc C. J. Amick}: Properties of steady Navier-Stokes solutions for certain unbounded channels and pipes. {\it Nonlinear Anal.} 2 (1978), no. 6, 689--720.

\bibitem{AP:1989SIAM} {\sc K. A. Ames and L. E. Payne}: Decay estimates in steady pipe flow. {\it SIAM J. Math. Anal.} 20 (1989), no. 4, 789--815.

\bibitem{BY:2020ANA} {\sc H. Beirao da Veiga and J. Yang}: Regularity criteria for Navier-Stokes equations with slip boundary conditions on non-flat boundaries via two velocity components. {\it Adv. Nonlinear Anal.}, 9 (2020), no. 1, 633--643,.

\bibitem{Bme1} {\sc M. Bogovski\u{i}}: Solution of the first boundary value problem for an equation of continuity of an incompressible medium, (Russian) {\it Dokl. Akad. Nauk SSSR} (248) 5 (1979), no. 3, 1037--1040.

\bibitem{Bme2} {\sc M. Bogovski\u{i}}: Solutions of some problems of vector analysis, associated with the operators div and grad, (Russian) {\it Theory of cubature formulas and the application of functional
analysis to problems of mathematical physics, pp. 5--40, 149, Trudy Sem. S. L. Soboleva, No. 1, 1980, Akad. Nauk SSSR Sibirsk. Otdel., Inst. Mat., Novosibirsk}, 1980.

\bibitem{CPZZ:2020ARMA} {\sc B. Carrillo, X. Pan, Q. S. Zhang and N. Zhao}: Decay and vanishing of some D-solutions of the Navier-Stokes equations. {\it Arch. Ration. Mech. Anal.} 237 (2020), no. 3, 1383--1419.

\bibitem{Chae:2014CMP} {\sc D. Chae}: Liouville-type theorems for the forced Euler equations and the Navier-Stokes equations. {\it Comm. Math. Phys.} 326 (2014), no. 1, 37--48.

\bibitem{CW:2016JDE} {\sc D. Chae and J. Wolf}: On Liouville type theorems for the steady Navier-Stokes equations in $\bR^3$. {\it J. Differential Equations} 261 (2016), no. 10, 5541--5560.

\bibitem{CQ:2010IUMJ} {\sc G.-Q. Chen and Z. Qian}: A study of the Navier-Stokes equations with the kinematic and Navier boundary conditions. {\it Indiana Univ. Math. J.}, 59 (2010), no. 2, 721--760.


\bibitem{Galdi:2011SPRINGER} {\sc G. P. Galdi}: {\it An introduction to the mathematical theory of the Navier-Stokes equations. Steady-state problems. Second edition.} Springer Monographs in Mathematics. Springer, New York, 2011.

\bibitem{GH:2011CMP} {\sc Y. Giga and H. Miura}: On vorticity directions near singularities for the Navier-Stokes fows with infinite energy, {\it Comm. Math. Phys.}, 303 (2011), 289--300.

\bibitem{GHM:2014CPDE} {\sc Y. Giga, P. Y. Hsu and Y. Maekawa}: A Liouville theorem for the planer Navier-Stokes equations with the no-slip boundary condition and its application to a geometric regularity criterion. {\it Comm. Partial Differential Equations} 39 (2014), no. 10, 1906--1935.


\bibitem{HW:1978SIAM} {\sc C. O. Horgan and L. T. Wheeler}: Spatial decay estimates for the Navier-Stokes equations with application to the problem of entry flow. {\it SIAM J. Appl. Math.} 35 (1978), no. 1, 97--116.

\bibitem{KNSS:2009ACTAMATH} {\sc G. Koch, N. Nadirashvili, G. A. Seregin and V. \v{S}ver\'{a}k}: Liouville theorems for the Navier-Stokes equations and applications. {\it Acta Math.} 203 (2009), no. 1, 83--105.


\bibitem{KO:1988RMS} {\sc  V.A. Kondrat'ev and O.A. Olenik}, Boundary value problems for a system in elasticity theory in unbounded domains. Korn inequalities. Russ. Math. Surv. 43 (1988), 65--119.

\bibitem{Konie:2006COLLMATH} {\sc P. Konieczny}: On a steady flow in a three-dimensional infinite pipe. {\it Colloq. Math.} 104 (2006), no. 1, 33--56.






\bibitem{Ladyzhenskaya:1959UMN} {\sc O. A. Lady\v{z}enskaya}: Investigation of the Navier-Stokes equation for stationary motion of an incompressible fluid. (Russian) {\it Uspehi Mat. Nauk} 14 (1959), no. 3, 75--97.

\bibitem{Ladyzhenskaya:1959SPD} {\sc O. A. Lady\v{z}enskaya}: Stationary motion of viscous incompressible fluids in pipes. {\it Soviet Physics. Dokl.} 124 (4) (1959), 68--70 (551--553 {\it Dokl. Akad. Nauk SSSR}).

\bibitem{Lady-Sol1980} {\sc O. Lady\v{z}enskaya and V. Solonnikov}: Determination of the solutions of boundary value problems for stationary Stokes and Navier-Stokes equations having an unbounded Dirichlet integral, Translated from Zapiski Nauchnykh Seminarov Leningradskogo Otdeleniya Matematicheskogo Instituta im. V. A. Steklova Akad. Nauk SSSR, Vol. 96, pp. 117-160, 1980.


\bibitem{Li-Schoen1984} {\sc P. Li and R. Schoen}: $L^p$ and mean value properties of subharmonic functions on Riemannian manifolds. {\it Acta Math.} 153(3--4), (1884), 279--301.

\bibitem{LPY:2022ARXIV} {\sc Z. Li, X. Pan and J. Yang}: On Leray's problem in an infinite-long pipe with the Navier-slip boundary condition. {\it arXiv}: 2204.10578.


\bibitem{PR:2017JFA} {\sc  D. Phan and S. S. Rodrigues.}: Gevrey regularity for Navier-Stokes equations under Lions boundary conditions. {\it J. Funct. Anal.}, 272 (2017), no. 7, 2865--2898.


\bibitem{KLW:2021ARXIV} {\sc K. Mikhail, W. Lyu and S. Weng}: On the existence of helical invariant solutions to steady Navier-Stokes equations. {\it arXiv}: 2102.13341.


\bibitem{MM:2009DIE}  {\sc M. Mitrea and S. Monniaux}: The nonlinear Hodge-Navier-Stokes equations in Lipschitz domains. {\it Differential Integral Equations,} 22 (2009), no. 3-4, 339--356.

\bibitem{Mucha:2003AAM} {\sc P. B. Mucha}: On Navier-Stokes equations with slip boundary conditions in an infinite pipe. {\it Acta Appl. Math.} 76 (2003), no. 1, 1--15.

\bibitem{Mucha:2003STUDMATH} {\sc P. B. Mucha}: Asymptotic behavior of a steady flow in a two-dimensional pipe. {\it Studia Math.} 158 (2003), no. 1, 39--58.

\bibitem{Navier} {\sc C. Navier}: Sur les lois du mouvement des fuides. {\it Mem. Acad. R. Sci. Inst. France,} 6, (1823), 389--440.

\bibitem{Pan:2020JMAA} {\sc X. Pan}: A Liouville theorem of Navier-Stokes equations with two periodic variables. {\it J. Math. Anal. Appl.} 485 (2020), no. 2, 123854, 7 pp.



\bibitem{Pileckas:2002MB} {\sc K. Piletskas}: On the asymptotic behavior of solutions of a stationary system of Navier-Stokes equations in a domain of layer type. (Russian) {\it Mat. Sb.} 193 (2002), no. 12, 69--104; translation in {\it Sb. Math.} 193 (2002), no. 11--12, 1801--1836.

\bibitem{Seregin:2016NON} {\sc G. Seregin}: Liouville type theorem for stationary Navier-Stokes equations. {\it Nonlinearity} 29 (2016), no. 8, 2191--2195.


\bibitem{Wang:2019JDE} {\sc W. Wang}:  Remarks on Liouville type theorems for the 3D steady axially symmetric Navier-Stokes equations. {\it J. Differential Equations} 266 (2019), no. 10, 6507--6524.


\bibitem{WX:2019ARXIV} {\sc Y. Wang and C. Xie}: Uniform structural stability of Hagen-Poiseuille flows in a pipe. {\it arXiv}: 1911.00749.


\bibitem{XX:2007CPAM} {\sc Y. Xiao and Z. Xin}: On the vanishing viscosity limit for the 3D Navier-Stokes equations with a slip boundary condition. {\it Comm. Pure Appl. Math.}, 60 (2007), no. 7, 1027--1055.


\bibitem{Zhang2021} {\sc Q. S. Zhang}: Bounded solutions to the axially symmetric Navier Stokes equation in a cusp region. {\it arXiv}:2106.08509 v1.

\end{thebibliography}
\end{document}